\documentclass[12pt]{article}
\usepackage{amsmath,latexsym,amssymb}
\usepackage{eucal}
\newtheorem{thm}{Theorem}[section]
\newtheorem{lem}[thm]{Lemma}
\newtheorem{prop}[thm]{Proposition}
\newtheorem{df}[thm]{Definition}

\newtheorem{cor}[thm]{Corollary}

\newcommand{\id}{\mathrm{id}}

\newcommand{\Coc}{C^1}

\DeclareMathOperator{\Aut}{Aut}

\DeclareMathOperator{\Ker}{Ker}

\DeclareMathOperator{\md}{mod}

\begin{document}

\title{Actions of discrete amenable groups into the normalizers of full groups of ergodic transformations}
\author{Toshihiko MASUDA \\
Faculty of Mathematics, Kyushu University \\
744, Motooka, Nishi-ku, Fukuoka 819-0395, Japan \\
E-mail: masuda@math.kyushu-u.ac.jp}
\date{}

\maketitle

\textit{
Dedicated to Professor Yasuyuki Kawahigashi on the occasion of his 60th birthday}

\begin{abstract}
 We apply Evans-Kishimoto's intertwining argument to the classification of 
actions of discrete amenable groups into the normalizer of a full group of an ergodic transformation.
Our proof does not depend on the types of ergodic transformations.
\end{abstract}

\section{Introduction}

The purpose of  this article is the  study actions of discrete amenable groups into the normalizer of a full group of an 
ergodic transformation on the Lebesgue space. 
The study of such objects has been motivated by 
the theory of operator algebras. In fact, many examples of von Neumann algebras arise from 
ergodic transformation through Krieger's construction. 

The study of automorphism groups of operator algebras  is one of the central subjects
for the theory of operator algebras,  and 
the classification of automorphisms and group actions 
has been developed since Connes' seminal works  \cite{Co-peri}, \cite{Con-auto}.
In particular  classification of actions of discrete amenable groups on injective factors
has been completed by many hands  \cite{J-act}, \cite{Ocn-act}, \cite{JT}, \cite{Su-Tak-RIMS},
\cite{Su-Tak-act}, \cite{KwST}, \cite{KtST}. These works heavily depend on 
the types of factors. However, we present the unified approach in \cite{M-unif-Crelle}
based on Evans-Kishimoto's method \cite{EvKi}, and 
gave proof that is independent of the types of factors.

There are corresponding results in ergodic theory.
The first result is due to  Connes and Krieger \cite{ConnKrieg}. 
They  developed the technique of use of ultraproduct to measure spaces and their transformation, and 
classified transformations (i.e., actions of $\mathbb{Z}$) in the normalizer of a full group of type II. 
Connes-Krieger's result has been generalized  
by \cite{Bezu-Golo-II} in the case of type II transformation and general discrete amenable groups,  
by  \cite{Bezu} in the case of type III$_\lambda$ transformation ($\lambda\ne 0$) and general discrete amenable groups
and finally by  \cite{Bezu-Golo-III0}
in the case of type III$_0$ transformation and general discrete amenable groups.
(See Theorem \ref{thm:main} below for the classification theorem.)
These results mentioned above depend on the types of transformations, and 
it is natural to  expect that our unified  approach \cite{M-unif-Crelle} is valid for classification of actions of discrete amenable groups 
into the normalizers of full groups on Lebesgue spaces. 
In fact, the answer is affirmative, and this is the main result of this article.
This classification result is very similar 
to that of the classification of actions of discrete amenable groups on injective factors. 
Indeed, classification result mentioned above  can be regarded as the 
classification of actions that fix Cartan subalgebras of 
Krieger factors. 

To apply the Evans-Kishimoto type intertwining argument,
we need  the characterization of full groups and their closures given by  
in  \cite{ConnKrieg} and \cite{Hamachi-JFA}. 
In the study of group actions on operator algebras, 
two classes of automorphisms play important roles, i.e.,
centrally trivial automorphisms and approximately inner automorphisms. 
In our case, full groups and their closures 
correspond to centrally trivial automorphism groups, and
approximately inner automorphism groups, respectively. 
Another main tool is the Rohlin type theorem,
Combining these results, we first show the cohomology vanishing theorem.
Then we obtain classification theorem by applying 
the Evans-Kishimoto type intertwining argument.

This paper is organized as follows. In \S \ref{sec:pre}, we collect basic facts which will be used in this paper,
and state the main results.
In \S \ref{Ultra}, we recall the ultraproduct construction of Connes-Krieger, and 
Ocneanu's Rohlin type theorem. 
In \S \ref{sec:2coho}, we show the second cohomology vanishing theorem. 
In \S \ref{sec:class}, we apply the Evans-Kishimoto type  intertwining argument \cite{EvKi}
 and classify actions of discrete amenable groups into the normalizer of a full group. 

\textbf{Acknowledgements.} 
The author is supported by JSPS KAKENHI Grant Number 16K05180, 
and 22K03341.

\section{Preliminaries}\label{sec:pre}

\subsection{Full groups of ergodic transformations and their normalizers}\label{subsec:topology}
In this subsection, we collect known facts on full groups of ergodic transformations and their normalizers,
which will be used in this article.

Let $(X,\mathcal{B},\mu)$ be a nonatomic Lebesgue  space with $\mu(X)=1$. 
(Throughout this article, we treat only nonatomic Lebesgue spaces.)
We denote by $\Aut(X,\mu)$ the set of all nonsingular transformations.
Fix an ergodic transformation $T\in \Aut(X,\mu)$. 
Let $[T]_*$ be  a set of all nonsingular bijection $R:A\rightarrow B$ for some $A, B\in \mathcal{B}$
such that $Rx\in \{T^nx\}_{n\in\mathbb{Z}}$, $x\in A$. 
Define a full group of $T$ by  $[T]:=[T]_*\cap \Aut(X,\mu)$, i.e., 
\[
[T]=\{R\in [T]_*\mid 
\mbox{the domain and the range of }R \mbox{ are both }X\}.
\]
We say $E, F\in \mathcal{B}$ are $T$-equivalent if there exists $R\in [T]_*$ whose domain is $E$ and 
range is $F$.  
A set $E\in \mathcal{B}$ is said to be 
$T$-infinite if there exists $F\subset E$ such that 
$\mu(E\backslash F)>0$  and  $F$ is $T$-equivalent to $E$. 
 A set $E\in \mathcal{B}$ is said to be $T$-finite if it is not $T$-infinite.

When $T$ is of type II, there exists the unique $T$-invariant measure $m$ on $X$ 
($m(X)<\infty$ when $T$ is of type II$_1$, and 
$m(X)=\infty$ when $T$ is of type II$_\infty$). In this case, 
the following two statements hold; (1) $E\in \mathcal{B}$ is $T$-finite if and only if $m(E)<\infty$, 
(2) $E,F\in \mathcal{B}$ are  $T$-equivalent if and only if $m(E)=m(F)$. 
When $T$ is of type II$_1$, we always assume $\mu$ is the unique $T$-invariant probability measure.

When $T$ is of type III, then any  $E\in \mathcal{B}$ with $\mu(E)>0$ is 
$T$-infinite. Hence if $E,F\in \mathcal{B}$ satisfy $\mu(E),\mu(F)>0$, then $E$ and $F$ 
are $T$-equivalent.

Let $N[T]\subset \Aut(X,\mu)$ be the normalizer of $[T]$. 
In the following, we use the notation $\hat{\alpha}(t)=\alpha t \alpha^{-1}$ for $t\in [T]$ and $\alpha\in N[T]$.

For $\alpha \in \Aut(X,\mu)$ and $\xi\in L^1(X,\mu)$, define $\alpha_\mu(\xi)\in L^1(X,\mu)$ by
\[
\alpha_\mu(\xi)(x):=\xi(\alpha^{-1}x)\frac{d (\mu\circ \alpha^{-1})}{d \mu}(x), \,\, \xi\in L^1(X,\mu).
\]
Then $\alpha_\mu$ is an isometry of $L^1(X,\mu)$, and $(\alpha\beta)_\mu=\alpha_\mu\beta_\mu$ holds for 
$\alpha,\beta\in \Aut(X,\mu)$. 

Let  $M(X,\mu)$ (resp. $M_1(X,\mu)$) be the set of 
complex-valued measures (resp. probability measures) which are absolutely continuous with respect to  $\mu$. 
For $\nu\in M(X,\mu)$, let $\|\nu\|=|\nu|(X)$, where $|\nu|$ is the total variation of $\nu$.
Then $M(X,\mu)$ is a Banach space with respect to the norm $\|\nu\|$.
For $\xi\in L^1(X,\mu)$, let $\nu_\xi(f)=\int_{X}\xi(x)f(x)d\mu(x)$.
Note that 
$L^1(X,\mu)$ and $M(X,\mu)$ are isomorphic as Banach spaces 
by $\xi\mapsto \nu_\xi$. 
Via this identification, $\alpha_\mu(\xi)$ corresponds to $\alpha(\nu_\xi)=\nu_\xi\circ \alpha^{-1}$.
In what follows, we freely use this identification, and we simply denote $\alpha_\mu(\xi)$ by $\alpha(\xi)$
for $\xi\in L^1(X,\mu)$. Thus $\xi(A)$, $A\in \mathcal{B}$, means $\nu_\xi(A)$.

Recall the topology of $N[T]$ introduced in \cite{HamachiOsikawaBook}. 
For $\alpha,\beta\in \Aut(X,\mu)$, 
$\{\alpha\ne \beta\}$ denotes the set $\{x\in X\mid \alpha x\ne\beta x\}$.
We say a sequence $\{\alpha_n\}_n\subset N[T]$ converges to $\beta\in N[T]$ weakly if 
$\lim\limits_{n\rightarrow \infty}\|\alpha_{n}(\xi)-\beta(\xi)\|=0$ for all $\xi\in M(X,\mu)$.
Define a metric $d_\mu$ by $d_\mu(\alpha,\beta):=\mu(\{\alpha \ne \beta \})$.
We say $\{\alpha_n\}_n\subset N[T]$ 
converges to $\beta\in N[T]$ uniformly if $\lim\limits_{n\rightarrow \infty}d_\mu(\alpha_n,\beta)=0$.
This definition does not depend on the choice of equivalence classes of $\mu\in M_1(X,\mu)$.
It is shown in \cite{HamaOsiSLN} that $[T]$ is a Polish group by $d_\mu$.

Now we gift a topology of $N[T]$ as follows. We say a sequence $\{\alpha_n\}_n\subset N[T]$ converges to  $\beta$ in $N[T]$ if 
$\{\alpha_n\}_n$ converges to $\beta$ weakly, and $\widehat{\alpha_n} (t)$ converges to $ \hat{\beta} (t) $
uniformly for all $t\in [T]$. (In fact, we only have to require convergence for $t\in \{T^n\}_{n\in \mathbb{Z}}$.)
This is the right topology for $N[T]$. In fact, this topology coincides with the 
$u$-topology for a Krieger factor $\mathcal{R}_{T}$ constructed from $(X,\mu, T)$.
So we also call this topology the $u$-topology.
It is shown that $N[T]$ is a Polish group in the $u$-topology \cite{HamachiOsikawaBook}.
Indeed, let $\{\xi_k\}_{k=1}^\infty\subset L^1(X,\mu)$ be a countable dense subset, and 
define a metric $d$ on $N[T]$ by 
\[
 d(\alpha,\beta):=
\sum_{k=1}^\infty \frac{1}{2^k}
\frac{\|\alpha(\xi_k)-\beta(\xi_k)\|}{1+\|\alpha(\xi_k)-\beta(\xi_k)\|}
+\sum_{k\in \mathbb{Z}}\frac{1}{2^{|k|}}
\frac{d_\mu(\hat{\alpha}(T^k),\hat{\beta}(T^k))}{1+d_\mu(\hat{\alpha}(T^k),\hat{\beta}(T^k))}.
\]
Then this $d$ makes $N[T]$ a Polish group, and 
the topology defined by $d$ is nothing but the $u$-topology on $N[T]$.

We collect elementary results which will be frequently used in what follows.
Since proof is easy,  we leave it to the readers
\begin{lem}\label{lem:inequality}
The following statements hold. \\
$(1)$ $d_\mu(\theta\alpha,\theta\beta)= d_\mu(\alpha,\beta)$, 
$d_\mu(\alpha\theta,\beta\theta)=  d_{\theta(\mu)}(\alpha,\beta)$, $\alpha,\beta,\theta\in N[T]$.
In particular we have $d_\mu(\alpha,\id)=d_\mu(\id,\alpha^{-1})=d_\mu(\alpha^{-1},\id)$, 
and $d_\mu(\hat{\alpha}(t),\hat{\alpha}(t'))=d_{\alpha^{-1}(\mu)}(t,t')$, $\alpha\in N[T]$, $t,t'\in[T]$. \\
$(2)$ $ d_{\nu_1}(\alpha,\beta)
\leq \|\nu_1-\nu_2\|+d_{\nu_2}(\alpha,\beta)$, $\nu_1,\nu_2\in M_1(X,\mu)$,
$\alpha,\beta\in N[T]$. \\
$(3)$ 
Let $\nu \in M_1(X,\mu)$, $A,B,C,D\in \mathcal{B}$. Then we have 
\begin{align*}
& \nu\left((A\cup B)\triangle (C\cup D)\right)\leq 
 \nu(A \triangle C)+
 \nu(B\triangle D), \\
 &\nu\left((A\cap B)\triangle (C\cap D)\right)\leq 
 \nu(A \triangle C)+
 \nu(B\triangle D).
\end{align*}
\end{lem}

Recall the definition of the fundamental homomorphism \cite{HamaOsiSLN}. 
Let $\tilde{X}:=X\times \mathbb{R}$, and $\mu_L$ be the Lebesgue measure on $\mathbb{R}$. 
For $R\in \Aut(X,\mu)$ and $t\in \mathbb{R}$, 
define $\tilde{R},F_t\in \Aut(\tilde{X},\mu \times \mu_L)$ by
\[
 \tilde{R}(x,u)=\left(Rx, u-\log\frac{d(\mu\circ R)}{d \mu}(x)\right),\,\,\, F_t(x,u)=(x,u+t).
\]
Let $(Y,\nu_Y)$ be the quotient space by $\tilde{T}$. Since $\tilde{T}$ and $F_t$ commute, 
we get the ergodic flow $(Y,\nu_Y,F_t)$, which is called the 
associated flow of $(X,T)$. 
Let 
\[
\Aut_F(Y,\nu_Y):=\{P\in \Aut(Y,\nu_Y)\mid PF_t=F_tP, t\in \mathbb{R}\}.
\]
When $R$ is in $N[T]$, $\tilde{R}$ induces $\md(R)\in \Aut_F(Y,\nu)$, which 
is called the fundamental homomorphism.
If we lift $R$ to an automorphism of a Krieger factor $\mathcal{R}_T$, $\md(R)$ is nothing but 
a Connes-Takesaki module for $R$ \cite{CT}.

In this article, we do not use the above definition of $\md(R)$ explicitly, and what we need is 
the fact $\Ker(\md)=\overline{[T]}$ (closure is taken in the $u$-topology) and the surjectivity of $\md$
\cite{HamaOsiSLN}, \cite{Hamachi-JFA}.

\subsection{Main results}\label{subsec:main}

\begin{df}
Let $G$ be a countable discrete group. \\
$(1)$  A map (or 1-cochain) $v:G\rightarrow [T]$ is said to be normalized if $v(e)=\id$.
We denote  the set of all normalized maps from $G$ into $[T]$ by $\Coc(G,[T])$. \\
$(2)$ A cocycle crossed action of $G$ into $N[T]$ is 
a pair of maps $\alpha: G\rightarrow N[T]$, and $c:G\times G\rightarrow [T]$ such that 
$\alpha_g\alpha_h= c(g,h)\alpha_{gh}$, $\alpha_e=\id$, $c(e,h)=c(g,e)=\id$. 
When $c(g,h)=\id$ for all $g,h\in G$, we say $\alpha$ is an action of $G$ into $N[T]$. \\
$(3)$ Let $(\alpha,c)$ be a cocycle crossed action of $G$ into $N[T]$, and $v\in C^1(G,[T])$.
A perturbed  crossed action $({}_v\alpha, {}_vc)$ of $(\alpha,c)$ by $v$ is defined by 
\[
{}_v\alpha_g:=v(g)\alpha_g,\,\,\, {}_vc(g,h)=v(g)\hat{\alpha_g}(v(h))c(g,h)v(gh)^{-1}.
\]
$(4)$ Let $\alpha$ be an action of $G$ into $N[T]$. We say 
a map $v\in \Coc(G,[T])$ 
is a 1-cocycle for $\alpha$ if 
$v$ satisfies the 1-cocycle identity $v(g)\widehat{\alpha_g}(v(h))=v(gh)$.
It is equivalent to that ${}_v\alpha$ is an action. \\ 
$(5)$ Let $\alpha$ and $\beta$ be actions of $G$ into $N[T]$. We say they are cocycle conjugate 
if there exist $\theta\in N[T]$ and 1-cocycle $v(\cdot )$
such that ${}_v\alpha_g=\theta\beta_g\theta^{-1}$ 
for all $g\in G$. If $\theta$ is chosen in $\overline{[T]}$, then we say they are strongly cocycle conjugate.
\end{df}

\noindent 
\textbf{Remark}
(1) Let $(\alpha,c)$ be a cocycle crossed action of $G$.
(Notion of a $p$-action is used in \cite{Bezu-Golo-III0}.)
By $(\alpha_g\alpha_h)\alpha_k=\alpha_g(\alpha_h\alpha_k)$, 
we can deduce the 2-cocycle identity $c(g,h)c(gh,k)=\widehat{\alpha_g}(c(h,k))c(g,hk)$. \\
(2) In many works, cocycle conjugacy is said to be outer conjugacy. 
In fact, we must distinguish these two notions for group actions on operator algebras, 
However,  in ergodic theory, we do not have to distinguish them.
(We have the canonical homomorphism $u\in [T]$ into the normalizer of a Krieger factor 
arising from $(X,\mu, T)$.)

\medskip 

At first, we show the following theorem.
\begin{thm}
Let $(\alpha,c)$ be a cocycle crossed action of a discrete amenable group into $N[T]$ with $\alpha_g\not\in [T]$,
$g\ne e$. 
Then $c(g,h)$ is a coboundary, that is, 
there exists  $v\in \Coc(G,[T])$
such that ${}_vc(g,h)=\id$, equivalently 
${}_v\alpha$ is a genuine action of $G$. 
If $c(g,h)$ is close to $\id$, then we can choose $v$ so that it is also close to $\id$.
\end{thm}
See below for a more precise statement. 

Let $N_\alpha:=\{g\in G\mid \alpha_g\in [T]\}$, which is a normal subgroup of $G$.
Our main result in this article is the following.
\begin{thm}\label{thm:main}
 Let $(X,\mu)$ be a Lebesgue  space with $\mu(X)=1$, $T$ an ergodic transformation on $(X,\mu)$.
Let $G$ be a countable discrete amenable group, and $\alpha$, $\beta$ actions of $G$ into $N[T]$. 
Then $\alpha$ and $\beta$ are strongly cocycle conjugate if and only if $N_\alpha=N_\beta$ and 
$\md(\alpha)=\md(\beta)$.  
\end{thm}
If $\alpha$ and $\beta$ are strongly cocycle conjugate, then it is obvious that 
$N_\alpha=N_\beta$ and $\md(\alpha_g)=\md(\beta_g)$. (Amenability of $G$ is unnecessary for this
implication.)  Thus the problem is to prove the converse implication, and a 
proof will be presented in 
subsequent sections. Here we only state the following corollary, which can be 
easily verified by Theorem \ref{thm:main}.

\begin{cor}
 Let $\alpha$ and $\beta$ be actions of $G$ into $N[T]$. Then $\alpha$ and $\beta$ are cocycle conjugate 
if and only if 
$N_\alpha$=$N_\beta$ and 
$\md(\alpha_g)=\theta\md(\beta_g)\theta^{-1}$ for some  $\theta\in \Aut_F(Y,\nu_Y)$.
\end{cor}
\textbf{Proof.} Since ``only if part'' is clear, we only have to prove ``if part''.
Suppose 
$N_\alpha$=$N_\beta$ and 
$\md(\alpha_g)=\theta\md(\beta_g)\theta^{-1}$ for some  $\theta\in \Aut_F(Y,\nu_Y)$.
By the surjectivity of $\md$ \cite{Hamachi-JFA}, we can take  $\sigma\in N[T]$ with 
$\md(\sigma)=\theta$. Then $\md(\alpha_g)=\md(\sigma\beta_g\sigma^{-1})$ holds, and hence
$\alpha_g$ and $\sigma\beta_g\sigma^{-1}$ are strongly cocycle conjugate 
by Theorem \ref{thm:main}. \hfill$\Box$

\section{Ultraproduct of a Lebesgue space and Rohlin type theorem}
\label{Ultra}
We recall ultraproduct construction in \cite{ConnKrieg}.

Let $\omega\in \beta \mathbb{N}$ be a free ultrafilter on $\mathbb{N}$.
For sequences $(A_n)_n, (B_n)_n\subset \mathcal{B}$, define an equivalence relation $(A_n)_n\sim (B_n)_n$ 
by $\lim\limits_{n\rightarrow \omega}\mu(A_n\triangle B_n)=0$. 
Let $\mathcal{B}^\omega:=\{(A_n)_n\subset \mathcal{B}\}\slash\!\!\sim$. 
This definition depends only on the equivalence class of $\mu$, 
and $\mathcal{B}^\omega$ is a boolean algebra. 

Any $\alpha \in N[T]$ induces a transformation $\alpha^\omega$ on $\mathcal{B}^\omega$
by $\alpha^\omega ((A_n)_n):=(\alpha (A_n))_n$. Let
\[
 \mathcal{B}_\omega:=\{\hat{A}\in \mathcal{B}^\omega: t^\omega \hat{A}=\hat{A}, t\in [T]\}.
\]
We denote by $\alpha_\omega$ the restriction of $\alpha^\omega$ on $\mathcal{B}_\omega$. 

Let $\hat{A}=(A_n)\in \mathcal{B}_\omega$. Then $\lim\limits_{n\rightarrow\omega}\chi_{A_n}$ 
exists in weak-$*$ topology on $L^\infty(X,\nu)$. By the ergodicity of $T$, 
this limit is in $\mathbb{C}$, and does not depend on the choice of representative $\hat{A}=(A_n)$. 
Thus we can define $\tau: \mathcal{B}_\omega\rightarrow \mathbb{C}$ by 
$\tau(A):=\lim\limits_{n\rightarrow \omega}\chi_{A_n}$. We can see $\tau\circ \alpha_\omega =\tau$ for $\alpha\in N[T]$.
  By \cite[Lemma 2.4]{ConnKrieg}, for $\alpha\in N[T]$, 
$\alpha_\omega=\id$ if and only if $\alpha\in [T]$. 
In fact, we have a stronger result. For $R\in N[T]$, if there exists $\hat{A}\in \mathcal{B}_\omega$ such that
$R_\omega \hat{B}=\hat{B}$ for any $\hat{B}\subset \hat{A}$, $\hat{B}\in \mathcal{B}_\omega$, 
then $R_\omega=\id$, and hence $R\in [T]$, \cite[Lemma 2.3]{ConnKrieg}. 
This means that $R_\omega$ is a free transformation if $R_\omega\ne \id$.

The main tool of this article is the following Rohlin type Theorem, essentially due to Ocneanu \cite{Ocn-act}.
(The following formulation is presented in \cite{M-unif-Crelle}.)

\begin{thm}\label{thm:Rohlin}
 Let $(\alpha,c)$ be a cocycle crossed  action of a discrete amenable group $G$ into $N[T]$ such that 
$\alpha_{g,\omega}\ne \id $ for all $g\ne e$. Let $K\Subset G$, $\varepsilon>0$, and $S$ be a
$(K,\varepsilon)$-invariant set. (The notation $K\Subset G$ means that $F$ is a finite subset of $G$.)
Then there exists a partition of unity $\{\hat{E}_s\}_{s\in S}\subset \mathcal{B}_\omega$ 
such that 
\begin{align*}
 (1)\,\,\, & \sum_{s\in S_g}\tau(\alpha_{g,\omega}\hat{E}_s\triangle \hat{E}_{gs})<5\varepsilon^{\frac{1}{2}},\,\, g\in K, \\
(2) \,\,\, & 
 \sum_{s\in S\backslash S_{g^{-1}}}\tau(\hat{E}_s)<3\varepsilon^{\frac{1}{2}}, 
\end{align*} 
where $S_g:=S\cap g^{-1}S$.
\end{thm}
Note that we have $gs\in S_{g^{-1}}=S\cap gS$ for  $s\in S_g$.

The proof of \cite{Ocn-act} is based on the following two facts, i.e., 
the freeness of actions on central sequence algebras, and ultraproduct technique. 
In our case, freeness holds as we remarked before Theorem \ref{thm:Rohlin}. 
Hence the proof of  \cite{Ocn-act} can be applied in our case by the suitable modification.

In what follows, we say $\alpha$ is an ultrafree action of $G$ if $\alpha_{g,\omega}\ne \id$ 
for any $g\in G$, $g\ne e$, 
to distinguish from the usual freeness of actions on Lebesgue spaces.

\begin{lem}\label{lem:choiceZs}
Let $A,B$ be finite sets,
$\{E_a\}_{a\in A}\subset \mathcal{B}_\omega$ a partition of $X$,  
and $\{P_{a,b}\}_{a \in A, b\in  B}\subset [T]$.
Choose representative $E_a=(E_a^n)_{n}$ such that $E_a^n \cap E_{a'}^n=\emptyset$ for $a\ne a'$, 
$\bigsqcup_{a\in A}E_a^n=X$. 
Then for any $\varepsilon>0$, $\Phi\Subset M_1(X,\mu)$, there exists $N\in \omega $, $\{Z_a^n\}_{a\in A}\subset \mathcal{B}$, 
$R^n_b \in [T]$, $n\in N $,  $b\in B$,
such that \\
$(1)$ $\nu(P_{a,b}^{-1}E_a^n \triangle E_a^n)<\varepsilon$, $n \in N$, $\nu \in \Phi$, \\
$(2)$ $Z_a^n \subset  E_a^n$, $P_{a,b}Z_a^n\subset E_a^n$, $n \in N$, \\
$(3)$ $\nu(E_a^n\backslash Z_a^n )<\varepsilon$, 
$\nu(E_a^n\backslash P_{a,b}Z_a^n )<\varepsilon$, 
$n \in N$, $\nu \in \Phi$, \\
$(4)$ $R_b^n x=P_{a,b}x$, $n\in N$, $x\in Z_a^n $. \\
\end{lem}
\textbf{Proof.} Since $P_{a,b}E_a=E_a$ by \cite[Lemma 2.4]{ConnKrieg}, there exists $N\in \omega  $ such that 
\[
 P_{a,b}(\nu)\left(\left(E_a^n \cup \bigcup_{b\in B}P_{a,b}^{-1}E_a^n \right) \backslash 
\left(E_a^n \cap \bigcap_{b\in B}P_{a,b}^{-1}E_a^n \right)\right)< 
\frac{\varepsilon}{2}
\]
for $n\in N$, $a\in A$, $b\in B$, $\nu\in \Phi$.

Let 
$Y_a^n:= E_a\cap \bigcap_{b\in B}P_{a,b}^{-1}E_a^n$.
Clearly we have  $Y_a^n, P_{a,b}Y_a^n\subset E_a^n$. Moreover
\[
\nu(P_{a,b}^{-1}E_a^n\triangle E_a^n)<\frac{\varepsilon}{2},\,\,
\nu(E_a^n\backslash Y_a^n)<\frac{\varepsilon}{2}, \,\, 
\nu(E_a^n\backslash P_{a,b}Y_a^n) =P_{a,b}(\nu)(P_{a,b}^{-1}E_a^n\backslash Y_a^n)
<\frac{\varepsilon}{2}\]
hold for $n \in N$, $\nu\in \Phi$. 
Let $Y^n:=\bigsqcup_{a\in A} Y_a^n$. 
Thus we can define $R^n_{0,b}\in [T]_*$ with $\mathrm{Dom}(R^n_{0,b})=Y^n$
by $R^n_{0,b}x=P_{a,b}x$, $x\in Y_a^n$.
If $X\backslash Y^n$ and $X\backslash R_{0,b}^n Y^n$ are $T$-equivalent, 
then we can extend $R_{0,b}^n$ to an element $R_b^n\in[T]$. 

At first, let us assume that $Y^n$ is T-finite. (Thus so is $R_{0,b}Y^n$.)
Such a case can happen if $T$ is of type II.
Then $X\backslash Y^n$ and $X\backslash R_{0,b}^n Y^n $ are $T$-equivalent.
Hence we can extend $R_{0,b}$ to $R_b\in [T]$. Set $Z_a^n:=Y_a^n$.
Then all the statements in the lemma are satisfied.

Next, let us assume that $Y^n$ is $T$-infinite. (Hence so is $R_{0,b}^nY^n$.) 
Take  $W_k\subset Y^n$, $k\in \mathbb{N}$,  such that $W_k\subset W_{k+1}$, $\bigcup_kW_k=Y^n$ and 
$Y^n \backslash W_k$ are $T$-infinite for all $k$. Set $Z_{a,k}^n:=Y_a^n\cap W_k$. 
Of course we have $Z_{a,k}^n\subset Z_{a,k+1}^n$, $\bigcup_{k}Z_{a,k}^n=Y_a^n$, 
$\bigsqcup_{a\in A}Z_{a,k}^n=Y^n \cap W_k=W_k$, 
and $Z_{a,k}^n, P_{a,b}Z_{a,k}^n \subset E_{a}^n$.
Thus $\{Z_{a,k}^n \}_{a\in A}$ satisfies the  condition (2).

Take sufficiently large $k$ such 
that 
\[
\nu(Y_a^n\backslash Z_{a,k}^n )<\frac{\varepsilon}{2}, \,\, 
\nu(P_{a,b}Y_a^n \backslash P_{a,b}Z_{a,k}^n )<\frac{\varepsilon}{2}  
\]
for $a\in A$, $b\in B$, $\nu \in \Phi$. 
Then it is clear that $\{Z_{a,k}^n \}$ satisfies the condition (3).
By the choice of $\{W_k\}$, $X\backslash \bigsqcup_{a\in A}Z_{a,k}^n \supset Y^n \backslash W_k$ 
and 
$X\backslash R_{0,b}^n\bigsqcup_{a\in A}Z_{a,k}^n \supset R_{0,b}^n(Y^n \backslash W_k)$.
It follows that 
$X\backslash \bigsqcup_{a\in A}Z_{a,k}^n$ and 
$X\backslash R_{0,b}^n\bigsqcup_{a\in A}Z_{a,k}^n $ are both $T$-infinite and hence are equivalent. 
Thus $Z_a^n:=Z_{a,k}^n$ satisfies all statements in the lemma. 
\hfill$\Box$

\medskip

Now we can combine Theorem \ref{thm:Rohlin} and Lemma \ref{lem:choiceZs} as follows.
\begin{prop}\label{prop:RohZschoice}
Let $G$ be a discrete amenable group, $(\alpha,c)$ an ultrafree cocycle crossed action of $G$ into $N[T]$.
Let $K \Subset G$ and $\varepsilon>0$ be given, and $S$ a $(K,\varepsilon)$-invariant set.
Let $B,C$ be finite set, 
$\{P_{s,b}\}_{s \in S, b\in  B}\subset [T]$, $\{\nu_s^c\}_{s\in S,c \in C}\Subset M_1(X,\mu)$. 
Then for any $\delta>0$, 
there exists a partition $\{E_s\}_{s\in S}\subset \mathcal{B}$ of $X$,
$E_s\supset Z_s$ and $R_{b}\in [T]$, $b\in B$, such that
\begin{align*}
  (1)\,\,\, & \sum_{s\in S_g}\nu_s^c (\alpha_gE_s\triangle E_{gs})<5\varepsilon^{\frac{1}{2}},\,\, g\in K,\, c\in C, \\
(2) \,\,\, &  \sum_{s\in S\backslash S_{g^{-1}}}\nu_s^c(E_s)<3\varepsilon^{\frac{1}{2}},  \,\, g\in K, \,c\in C,  \\ 
(3)\,\,\,  & \nu_s^c(P_{s,b}^{-1}E_s^n \triangle E_s^n)<\delta, \,\, s \in S, \,b\in B, \, c\in C,\\
 (4)\,\,\,& P_{s,b}Z_a\subset E_s,\,\, s\in S, \, b\in B,  \\
 (5)\,\,\, &\nu_s^c(E_s\backslash Z_s)<\delta, \nu_s^c(E_s\backslash P_{s,b}Z_s)<\delta, 
\,\, s\in S, \,b\in B, \, c\in C,\\
 (6)\,\,\,& R_b x=P_{s,b}x,\,\, s\in S, \, b\in B, \, x\in Z_s. 
\end{align*}
\end{prop}
\textbf{Proof.} Let $\{\hat{E}_s\}_{s\in S}\subset \mathcal{B}_\omega$ be a Rohlin partition as in 
Theorem \ref{thm:Rohlin}. 
Since $\tau(\hat{A})=\lim\limits_{n\rightarrow \omega}\chi_{A_n}$ for $\hat{A}=(A_n)_n\in \mathcal{B}_\omega$,
$\tau(\hat{A})=\lim\limits_{n\rightarrow \omega}\nu(A^n)$ for any $\nu\in M_1(X,\mu)$.
Choose representative $\hat{E}_s=(E^n_s)_n$ such that 
$E_s^n\cap E_{s'}^n=\emptyset$, $\bigsqcup_{s\in S}E_s^n=X$. 
By Theorem \ref{thm:Rohlin}, 
\begin{align*}
 (1)\,\,\, & \lim\limits_{n\rightarrow \omega}\sum_{s\in S_g}\nu_s^c(\alpha_gE_s^n\triangle 
 E_{gs}^n)<5\varepsilon^{\frac{1}{2}},\,\, g\in K, \\
(2) \,\,\, & \lim_{n\rightarrow \omega}\sum_{s\in S\backslash S_{g^{-1}}}\nu_s^c(E_s^n)<3\varepsilon^{\frac{1}{2}},
 \,\, g\in K 
\end{align*} 
holds  for any $\{\nu_s^c\}_{s\in S,c\in C}\subset M_1(X,\mu)$. 
Thus there exists $N_1\in \omega$ such that
\begin{align*}
 & \sum_{s\in S_g}\nu_s^c(\alpha_gE_s^n\triangle E_{gs}^n)<5\varepsilon^{\frac{1}{2}},\,\, g\in K,\, c\in C, \\
 & \sum_{s\in S\backslash S_{g^{-1}}}\nu_s^c(E_s^n)<3\varepsilon^{\frac{1}{2}}, \,\, g\in K,\, c\in C
\end{align*} 
for all $n\in N_1$. By Lemma \ref{lem:choiceZs}, there exists $N_2\in \omega$,
$Z_s^n\subset E_s^n$ and $R_b^n\in [T]$, ($n\in N_2$),  
 such that  
\begin{align*}
 & \nu_s^c(P_{s,b}^{-1}E_s^n \triangle E_s^n)<\delta, \,\, s\in S,  \, b\in B, \\
  & P_{s,b}Z_s^n\subset E_s^n,\,\, s\in S, b\in B,  \\
  &\nu_s^c(E_s^n\backslash Z_s^n)<\delta, \nu_s^c(E_s^n\backslash P_{s,b}Z_s^n)<\delta, \,\, s\in S,\, b\in B,\, c\in C, \\
 & R_b^n x=P_{s,b}x,\,\, s\in S, \, b\in B,\, x\in Z_s^n
\end{align*}
for any $n\in N_2$.
Fix $n\in N_1\cap N_2$, and set $E_s:=E_s^n$, $Z_s:=Z_s^n$, $R_b:=R_b^n$.
Then these $E_s, Z_s, R_b$ are desired objects. \hfill$\Box$

\section{Cohomology vanishing}\label{sec:2coho}

At first, we show the following second cohomology vanishing result, which is shown 
in \cite[Theorem 1.3]{Bezu-Golo-III0}. We present the proof for readers' convenience.

\begin{thm}\label{2cohoinfinite}
 Let $T$ be a transformation of type II$_\infty$ or type III, and $(\gamma,c)$ a cocycle crossed action 
of a discrete group $G$ into $N[T]$. Then $c(g,h)$ is a coboundary, i.e.,
there exists $u\in \Coc(G,[T])$
such that ${}_uc(g,h)=\id$.  
\end{thm}
\textbf{Proof.} Since $T$ is of type II$_\infty$, or type III, there exists a partition 
$\{E_h\}_{h\in G}$ of $X$ such that each $E_h$ is $T$-infinite.
Let $\{f_{g,h}\}_{g,h\in G}\subset [T]$ be an array for $\{E_g\}_{g\in G}$.
Take $v_g^0 \in [T]_*$ with $\mathrm{Dom}(v_g^0)=\gamma_gE_e$ and $\mathrm{Ran}(v_g^0)=E_e$.
Define $v(g)\in [T]$ by $ f_{h,e}v_g^0\gamma_g(f_{e,h})$ on $\gamma_gE_h$.
Then we have ${_v}\gamma_g:E_h\rightarrow E_h$ 
and 
$\widehat{{}_v\gamma_g}(f_{h,k})=f_{h,k}$ for any $g,h,k\in G$. 
Replacing $(\gamma,c)$ with $({}_v\gamma,{}_vc)$, we may
assume $\gamma_gE_k=E_k$ and $\widehat{\gamma_g}(f_{h,k})=f_{h,k}$. 
Since $\gamma_g\gamma_h=c(g,h)\gamma_{gh}$, we also have  $c(g,h)E_k=E_k$ and 
$\widehat{c(g,h)}(f_{k,l})=f_{k,l}$.

Next define $u(g)\in [T]$ by $u(g)=c(g,l)^{-1}f_{gl, l}$ on $E_{l}$.
Note $u(g)$ sends $E_{l}$ to $E_{gl}$, hence so does ${}_u\gamma_g$.
Hence for $x\in E_l$,
\begin{align*}
 {}_u\gamma_g\, {}_u\gamma_h x &=u(g)\gamma_g c(h,l)^{-1}  f_{hl,l}\gamma_h x 
=u(g)\gamma_g c(h,l)^{-1} \gamma_{h} f_{hl,l}  x \\
&=u(g)\widehat{\gamma_g}(c(h,l))^{-1} \gamma_g\gamma_{h} f_{hl,l}  x 
=u(g) \widehat{\gamma_g}(c(h,l))^{-1} c(g,h)\gamma_{gh} f_{hl,l}  x \\
&=c(g,hl)^{-1}f_{ghl,hl} \widehat{\gamma_g}(c(h,l))^{-1} c(g,h)\gamma_{gh} f_{hl,l}  x \\
&=c(g,hl)^{-1} \widehat{\gamma_g}(c(h,l))^{-1} c(g,h)\gamma_{gh}f_{ghl,hl} f_{hl,l}  x \\
 &=c(gh,l)^{-1}f_{ghl,l}\gamma_{gh}x =u(gh)\gamma_{gh} x. 
\end{align*}
This implies that  ${}_u\gamma$ is an action, and ${}_uc(g,h)=u(g)\widehat{\gamma_g}(u(h))c(g,h)u(gh)^{-1}=\id$ holds.
 \hfill $\Box$
 
\medskip

In the Theorem \ref{2cohoinfinite}, we have no estimation on the choice of $u(g)$, even if 
$c(g,h)$ is close to $\id$. The rest of this section is devoted to solving this problem.
From now on, we always assume that $G$ is a discrete amenable group.

\medskip

For all $g\in G$ and  $S\Subset G$, 
fix a bijection $l(g): S\rightarrow S$ such that $l(g)s=gs$ if $gs\in S$. 

\begin{lem}\label{lem:appro2cochovanish1}
Let $(\gamma,c)$ be an ultrafree  cocycle crossed action of $G$. 
For any $\varepsilon>0$, $K\Subset G$, $\mu\in \Phi\Subset M_1(X,\mu)$, 
there exists $w\in \Coc(G,[T])$
such that
\[
d_\nu({}_wc(g,h),\id)<\varepsilon,\,\, g,h \in K,\, \nu \in \Phi.
\]
Moreover for given $\varepsilon>0$, 
$e\in K\Subset G$, 
there exist $\delta>0$ and $S\Subset G$, which depends only on $K$ and $\varepsilon>0$, 
such that if 
\[
 \|c(g,h)(\xi)-\xi\|<\delta,\,\, 
d_\nu(\widehat{c(g,h)}(t),t)<\delta,\,\, g,h\in S,\, t\in \Lambda ,\, \xi,\nu\in \Phi
\]
for  some cocycle crossed action $(\gamma,c)$, $\Lambda \Subset [T]$ and $\Phi\Subset M_1(X,\mu)$, 
then we can choose $w\in \Coc(G,[T])$ so that it  further satisfies
\[
 \|w(g)(\xi)-\xi\|<\varepsilon,\,\, d_\nu(\widehat{w(g)}(t),t)<\varepsilon,\,\, g\in K, \,\xi,\nu\in \Phi, \,t\in \Lambda .
\]
\end{lem}
\textbf{Proof.}
Choose $\varepsilon'>0$ with $11\sqrt{\varepsilon'}<\varepsilon$, and let $S'\subset G $ be 
a $(K\cup K^2,\varepsilon')$-invariant set, and $S=S'\cup K$. 
Choose $\delta$ such that 
$5\delta|S|+11\sqrt{\varepsilon'}<\varepsilon$. 

By  applying Proposition \ref{prop:RohZschoice},
we can take Rohlin partition $\{E_s\}_{s\in S'}\subset \mathcal{B}$, $Z_s\subset E_s$, $w(g)\in [T]$, $g\in K$, such that 
\begin{align*}
(1)\,\,\, & E_{l(g)s}\supset c(g,s)^{-1}Z_{l(g)s},\,\, g\in K\cup K^2, \,s\in S', \\
(2) \,\,\, & 
 \nu(E_s\backslash Z_s)<\delta,\,\, \nu(E_{l(g)s}\backslash  c(g,s)^{-1}Z_{l(g)s})<\delta, \,\, 
    g\in K\cup K^2, \,s\in S',\, \nu \in \Phi, \\
(3)\,\,\, &  \nu\left(c(gh,k)^{-1}c(g,h)^{-1} \widehat{\gamma_g}(c(h,k))(E_{ghk}\backslash  Z_{ghk})\right) <\delta, 
\,\,\, g,h\in K,\, k\in S'_{gh}\cap S'_h, \,\nu\in \Phi, \\
(4)\,\,\, &  \nu\left(c(gh,k)^{-1}c(g,h)^{-1}\gamma_g\left(E_{hk}\backslash  Z_{hk}\right)\right) <\delta, 
\,\, g,h\in K, \,k\in  S'_h,\,\nu \in \Phi, \\
(5)\,\,\,& \nu \left(E_{ghk}\triangle c(gh,k)^{-1}c(g,h)^{-1} \widehat{\gamma_g}(c(h,k))E_{ghk}\right) <\delta,
\,\, g,h\in K, \,k\in  S'_h,\,\nu\in \Phi, \\
(6)\,\,\,&  \nu(E_{ghk}\triangle c(gh,k)^{-1}c(g,h)^{-1}E_{ghk})<\delta, \,\, g,h\in K, \,k\in S'_{gh}\cap S'_h, \, \nu \in \Phi, \\
(7)\,\,\, &\sum_{k\in S'_{gh}\cap S'_h}
\nu\left(c(gh,k)^{-1}c(g,h)^{-1}\left(E_{ghk}\triangle \gamma_gE_{hk}\right)\right) <5\sqrt{\varepsilon'}, 
 \,\, g,h\in K,\,\nu\in \Phi, \\ 
(8)\,\,\,& \sum_{k\in S'\backslash S'_{(gh)^{-1}}}\nu( E_s)<3\sqrt{\varepsilon'}\,\, g\in K\cup K^2,\, \nu \in \Phi, \\
(9)\,\,\, & w(g)x=c(g,s)^{-1}x, \,\, x\in Z_{l(g)s},\, g\in K, \,s\in S'.
\end{align*}
Here we applied Proposition \ref{prop:RohZschoice} for 
\begin{align*}
 B&=\{c(g,s)^{-1}\mid g\in K, s\in S'\}
\cup \{c(gh,k)^{-1}c(g,h)^{-1}\mid g,h\in K, k\in S'\} \\
&\cup \{c(gh,k)^{-1}c(g,h)^{-1}\widehat{\gamma_g}(c(h,k))\mid g,h\in K, k\in S'\},
\end{align*}
and 
\begin{align*}
 C&=
\Phi \cup \left\{
\nu\left(c(gh,k)^{-1}c(g,h)^{-1} \widehat{\gamma_g}(c(h,k))\,\,\cdot\,\, \right)\mid \nu\in \Phi, g,h\in K, k\in S' \right\} \\
&\cup 
\left\{
\nu\left(c(gh,k)^{-1}c(g,h)^{-1} \gamma_g\,\,\cdot\,\, \right)\mid \nu\in \Phi, g,h\in K, k\in S' \right\} \\
&\cup 
\left\{
\nu\left(c(gh,k)^{-1}c(g,h)^{-1} \,\, \cdot\,\, \right)\mid \nu\in \Phi, g,h\in K, k\in S' \right\}.
\end{align*}

We define $w(g)=\id$ if $g\not\in G$.

Let 
 \begin{align*}
W^0_{g,h,k} = c(gh,k)^{-1}Z_{ghk}\cap c(gh,k)^{-1}c(g,h)^{-1}\gamma_gZ_{hk}\cap 
c(gh,k)^{-1}c(g,h)^{-1} \widehat{\gamma_g}(c(h,k))Z_{ghk}
 \end{align*}
for $k\in S'_{g,h}\cap S'_h $ and 
\[
 W_{g,h}= \bigcup_{k\in S'_{gh}\cap S'_{h}} W^0_{g,h,k}.
\]
We can verify $w(g)\hat{\gamma_g}(w(h))c(g,h)w(gh)^{-1}=\id$ on $W_{g,h}$ as follows.
Take $x\in W^0_{g,h,k}$. Since $x\in c(gh,k)^{-1}Z_{ghk}$, we have
$w(gh)^{-1}x=c(gh,k)x$. 
Thus we have $\gamma_g^{-1}c(g,h)w(gh)^{-1}x=\gamma_g^{-1}c(g,h)c(gh,k)x$.
Since $x\in c(gh,k)^{-1}c(g,h)^{-1}\gamma_gZ_{hk}$, 
$\gamma_g^{-1}c(g,h)c(gh,k)x\in Z_{hk}$ holds.
Hence we have
\[
 w(h)\gamma_g^{-1}c(g,h)c(gh,k)x= c(h,k)^{-1}\gamma_g^{-1}c(g,h)c(gh,k)x.
\]
Since $x\in c(gh,k)^{-1}c(g,h)^{-1} \widehat{\gamma_g}(c(h,k))Z_{ghk}$,
 \[
\widehat{  \gamma_g}(w(h))c(g,h)w(gh)^{-1}x=\widehat{\gamma_g}(c(h,k))^{-1}c(g,h)c(gh,k)x\in Z_{ghk}
 \] 
holds, and hence we have
\[
 w(g)\gamma_gc(h,k)^{-1}\gamma_g^{-1}c(g,h)c(gh,k)x=
 c(g,hk)^{-1}\gamma_gc(h,k)^{-1}\gamma_g^{-1}c(g,h)c(gh,k)x=x
\]
by the 2-cocycle identity. This shows 
${}_wc(g,h)=w(g)\widehat{\gamma_g}(w(h))c(g,h)w(gh)^{-1}=\id$ on $W_{g,h}$.
Thus we have $\{{}_wc(g,h)\ne \id\}\subset X\backslash W_{g,h}$.

We will show $\nu(X\backslash  W_{g,h})<\varepsilon$ for $\nu\in \Phi$. 
By (2), we have
\[
 \nu(E_{ghk}\backslash c(gh,k)^{-1}Z_{ghk})<\delta,\,\,\,
\nu\in \Phi, \, g,h\in K\, k\in S'_{gh}.
\]

For $g,h\in K$, $k\in S_{gh}'\cap S_h'$, $\nu\in \Phi$, we have
\begin{align*}
\lefteqn{ \nu(E_{ghk}\triangle c(gh,k)^{-1}c(g,h)^{-1}\gamma_gZ_{hk})} \\
&\leq 
 \nu(E_{ghk}\triangle c(gh,k)^{-1}c(g,h)^{-1}E_{ghk})+
\nu\left(c(gh,k)^{-1}c(g,h)^{-1}\left(E_{ghk}\triangle \gamma_gZ_{hk}\right)\right) \\
&\leq \delta+
\nu\left(c(gh,k)^{-1}c(g,h)^{-1}\left(E_{ghk}\triangle \gamma_gZ_{hk}\right)\right) \,\,\, (\mbox{by }(6)) \\
&\leq \delta+
\nu\left(c(gh,k)^{-1}c(g,h)^{-1}\left(E_{ghk}\triangle \gamma_gE_{hk}\right)\right)+
\nu\left(c(gh,k)^{-1}c(g,h)^{-1}\gamma_g\left(E_{hk}\triangle Z_{hk}\right)\right) \\
&\leq 2\delta+
\nu\left(c(gh,k)^{-1}c(g,h)^{-1}\left(E_{ghk}\triangle \gamma_gE_{hk}\right)\right) \,\,\, (\mbox{by }(4))
\end{align*}
and \begin{align*}
\lefteqn{\nu \left(E_{ghk}\triangle c(gh,k)^{-1}c(g,h)^{-1} \widehat{\gamma_g} (c(h,k))Z_{ghk}\right)} \\
&\leq
 \nu \left(E_{ghk}\triangle c(gh,k)^{-1}c(g,h)^{-1} \widehat{\gamma_g}(c(h,k))E_{ghk}\right)  \\
&\hspace*{12pt}+
\nu \left(c(gh,k)^{-1}c(g,h)^{-1} \widehat{\gamma_g}(c(h,k))(E_{ghk}\triangle Z_{ghk})\right) \\
&<2\delta \,\,\, (\mbox{by (5) and (3)}).
\end{align*}

Thus
\begin{align*}
 \nu(E_{ghk}\triangle  W_{g,h,k}^0)&\leq \nu(E_{ghk}\backslash c(gh,k)^{-1}Z_{ghk}) 
 +\nu(E_{ghk}\triangle c(gh,k)^{-1}c(g,h)^{-1}\gamma_gZ_{hk})\\
&\hspace*{12pt} +\nu \left(E_{ghk}\triangle c(gh,k)^{-1}c(g,h)^{-1} \gamma_gc(h,k)\gamma_g^{-1}Z_{ghk}\right) \\
&<5\delta+ 
\nu\left(c(gh,k)^{-1}c(g,h)^{-1}\left(E_{ghk}\triangle \gamma_gE_{hk}\right)\right)
\end{align*}
follows. Then 
\begin{align*}
\nu\left( \left(\bigcup_{k\in S'_{gh}\cap S'_h}E_{ghk}\right)\triangle W_{g,h}\right) 
&\leq 
\sum_{k\in S'_{gh}\cap S'_h}
\nu\left(  E_{ghk} \triangle W_{g,h,k}^0\right) \\
&< 5\delta|S|+ 
\sum_{k\in S'_{gh}\cap S'_h}
\nu\left(c(gh,k)^{-1}c(g,h)^{-1}\left(E_{ghk}\triangle \gamma_gE_{hk}\right)\right) \\
&<5\delta|S|+5\sqrt{\varepsilon'} \,\,\, (\mbox{by (7)})
\end{align*}
holds. 

Finally, we have 
\begin{align*}
 \nu(X\backslash W_{g,h})&\leq 
\nu\left(X\backslash \bigcup_{k\in S'_{gh}\cap S'_h}E_{ghk}\right)
+\nu\left( \left(\bigcup_{k\in S'_{gh}\cap S'_h}E_{ghk}\right)\triangle W_{g,h}\right) \\
&=\nu\left(\bigcup_{k\in S\backslash  \left(S'_{g^{-1}}\cap S'_{(gh)^{-1}}\right)}E_{k} \right)
+\nu\left( \left(\bigcup_{k\in S'_{gh}\cap S'_h}E_{ghk}\right)\triangle W_{g,h}\right) \\
&<\sum_{k\in S'\backslash S'_{g^{-1}}}\nu(E_k)+\sum_{k\in S'\backslash S'_{(gh)^{-1}}}\nu( E_s)+
5\delta|S|+5\sqrt{\varepsilon'} \\
&<5\delta|S|+11\sqrt{\varepsilon'} \,\,\, (\mbox{by (8)}) \\
&<\varepsilon.
\end{align*}
(Note $ghk\in S'_{g^{-1}}\cap S'_{(gh)^{-1}}$ for $k\in S_h'\cap S_{gh}'$.)
This implies 
$d_\nu({}_wc(g,h),\id)<\varepsilon$ for $g,h\in K$, $\nu\in \Phi$. 

Assume 
\[
 \|c(g,h)(\xi)-\xi\|<\delta,\,\, 
d_\nu(\widehat{c(g,h)}(t),(t))<\delta,\,\, g,h\in S,\, t\in \Lambda ,\, \xi,\nu\in \Phi.
\]
We show 
\[
 \|w(g)(\xi)-\xi\|<\varepsilon,\,\, d_\nu(\widehat{w(g)}(t),t)<\varepsilon,\,\, 
g\in K, \,\xi,\nu\in \Phi, \,t\in \Lambda .
\]

Let $Z:=\bigsqcup_{s\in S'}Z_s$. By the definition of $w(g)$, $w(g)Z=\bigsqcup_{s\in S'}c(g,s)^{-1}Z_s$ holds.
Then we have 
\begin{align*}
\lefteqn{ \|w(g)(\xi)-\xi\|} \\
 &=\int_{X}|w(g)(\xi)(x)-\xi(x)|d\mu(x)  \\
&=\int_{w(g) Z}|w(g)(\xi)(x)-\xi(x)|d\mu(x)  +\int_{X\backslash w(g)Z}|w(g)(\xi)(x)-\xi(x)|d\mu(x)  \\
&=
\sum_{s\in S'}\int_{c(g,s)^{-1}Z_{l(g)s}}|w(g)^{-1}(\xi)(x)-\xi(x)|d\mu(x)
+\int_{X\backslash w(g)Z}|w(g)(\xi)(x)-\xi(x)|d\mu(x). 
\end{align*}

Note 
\[
 w(g)(\xi)(x)=\xi(w(g)^{-1}x)\frac{d(\mu \circ w(g)^{-1})}{d\mu}(x)
=\xi(c(g,s)x)\frac{d(\mu\circ c(g,s))}{d\mu}(x)=c(g,s)^{-1}(\xi)(x)
\]
for $x\in c(g,s)^{-1}Z_{l(g)s}$, when we regard $\xi$ as an element of $L^1(X,\mu)$. Thus 
the first term is estimated as follows;
\begin{align*}
\lefteqn{\sum_{s\in S'}\int_{c(g,s)^{-1}Z_{l(g)s}}|w(g)^{-1}(\xi)(x)-\xi(x)|d\mu(x)} \\
&=
\sum_{s\in S'}\int_{c(g,s)^{-1}Z_{l(g)s}}|c(g,s)^{-1}(\xi)(x)-\xi(x)|d\mu(x) \\
&\leq \sum_{s\in S'}\|c(g,s)^{-1}(\xi)-\xi\| <\delta|S|.
\end{align*}

To estimate the second term, one should note 
\[
 \xi(X\backslash Z)=\sum_{s\in S'}\xi(E_s\backslash Z_s)<\delta|S'|,\,\,\,
 \xi(X\backslash w(g)Z)=\sum_{s\in S'}\xi(E_s\backslash c(g,s)^{-1}Z_s)<\delta|S'|
\] 
by (2). Hence
\begin{align*}
\lefteqn{\int_{X\backslash w(g)Z} |w(g)(\xi)(x)-\xi(x)|d\mu(x) } \\
&\leq \int_{X\backslash w(g)Z}w(g)(\xi)(x) d\mu(x)+
\int_{X\backslash w(g)Z}\xi(x) d\mu(x) \\
&=\int_{X\backslash w(g)Z}\xi(w(g)^{-1}x)\frac{d(\mu\circ  w(g)^{-1})}{d\mu} d\mu(x)+\xi(X\backslash w(g)Z) \\
&= \int_{X\backslash Z}\xi(x)d\mu(x)+\xi(X\backslash w(g)Z) 
= \xi(X\backslash Z)+\xi(X\backslash w(g)Z) <2|S'|\delta,
 \end{align*}
and we obtain $\|w(g)(\xi)-\xi\|<3\delta |S'|<\varepsilon$.

We next show 
\[
 d_\nu(\widehat{w(g)}(t),t)<\varepsilon,\,\, g\in G, \,\nu\in \Phi, \,t\in \Lambda.
\]
By the assumption 
\[
 \|c(g,s)(\nu)-\nu\|<\delta,\,\, 
d_\nu(\widehat{c(g,s)}(t),t)<\delta,\,\, t\in \Lambda ,\, g,s\in S, \,\nu \in \Phi, 
\]
\begin{align*}
 d_\nu(\widehat{c(g,s)^{-1}}(t), t)&=d_{c(g,s)(\nu)}(t, \widehat{c(g,s)}(t)) \\
&\leq \|c(g,s)(\nu)-\nu\|+d_{\nu}(\widehat{c(g,s)}(t),t) 
< 2\delta
\end{align*}
holds.
We can further assume 
\[
 \nu\left(E_{l(g)s}\triangle  c(g,s)^{-1}t^{-1}E_{l(g)s}\right)<\delta, \,\,
 \nu\left(c(g,s)^{-1}t^{-1}\left(E_{l(g)s}\triangle Z_{l(g)s}\right)\right) < \delta\]
for $t\in \Lambda $, $s\in S'$, $g\in K$ in  the choice of $Z_s$ and $E_s$.

Let $B_{g,s,t}:=\{\widehat{c(g,s)^{-1}}(t)=t\}$. Then we have $\nu(X\backslash B_{g,s,t})<2\delta$, $\nu\in \Phi$.
We can see $\widehat{w(g)}(t) =\widehat{c(g,s)^{-1}}(t)$ on 
$c(g,s)^{-1}Z_{l(g)s}\cap c(g,s)t^{-1}Z_{l(g)s}$ as above.
Thus  $\widehat{w(g)}(t)=t$ holds on
\[
 \bigcup_{s\in S'}c(g,s)^{-1}Z_{l(g)s}\cap c(g,s)t^{-1}Z_{l(g)s}\cap B_{g,s,t}.
\]
We will show
\[
 \nu\left(X\backslash  \bigcup_{s\in S'}c(g,s)^{-1}Z_{l(g)s}\cap c(g,s)t^{-1}Z_{l(g)s}\cap B_{g,s,t}
 \right)<\varepsilon.
\]
At first, we have
\begin{align*}
\lefteqn{\nu\left(E_{l(g)s}\triangle  \left(c(g,s)^{-1}Z_{l(g)s}\cap c(g,s)^{-1}t^{-1}Z_{l(g)s}\right)\right)} \\
&\leq 
\nu(E_{l(g)s}\triangle c(g,s)^{-1}Z_{l(g)s})+
\nu(E_{l(g)s}\triangle c(g,s)^{-1}t^{-1}Z_{l(g)s}) \\
&<\delta +
\nu(E_{l(g)s}\triangle  \cap c(g,s)^{-1}t^{-1}E_{l(g)s})+
\nu\left(c(g,s)^{-1}t^{-1}\left(E_{l(g)s}\triangle Z_{l(g)s}\right)\right)\,\, \mbox{(by (2))} \\
&<3\delta.
\end{align*}

Thus 
\begin{align*}
\lefteqn{ \nu \left(X\backslash \bigcup_{s\in S'}c(g,s)^{-1}Z_{l(g)s}\cap c(g,s)^{-1}t^{-1}Z_{l(g)s}\cap B_{g,s,t}\right) } \\
&\leq 
 \nu \left(X\backslash \bigcup_{s\in S'}c(g,s)^{-1}Z_{l(g)s}\cap c(g,s)^{-1}t^{-1}Z_{l(g)s}\right) 
+ \nu\left(X\backslash \bigcup_{s\in S'}B_{g,s,t}\right) \\
&\leq 
\sum_{s\in S'}
 \nu\left(E_{l(g)s}\triangle  \left(c(g,s)^{-1}Z_{l(g)s}\cap c(g,s)^{-1}t^{-1}Z_{l(g)s}\right)\right) + 
\sum_{s\in S'}\nu\left(X\backslash B_{g,s,t}\right) \\
&<5|S'|\delta<\varepsilon
\end{align*}
holds, and we obtain $d_\nu(\widehat{w(g)}(t),t)<\varepsilon$ for $g\in K$, $\nu\in \Phi$, $t\in \Lambda$. \hfill$\Box$

\begin{lem}\label{lem:near1cocyclevanish}
 For any $e\in K\Subset  G $ and  $\varepsilon>0$, there exist $S\Subset G$ and $\delta>0$ 
satisfying the following property;
for any $\mu\in \Phi\Subset M_1(X,\mu)$, an ultrafree cocycle crossed action $(\gamma, c)$ of $G$, and  
$u\in \Coc(G,[T])$   with 
\[
d_\nu(\widehat{\gamma_g}(u(s))^{-1}u(g)^{-1}u(gs), \id )<\delta, \,\, g,s\in S, \,\nu\in \Phi, 
\]
there exists $w\in [T]$ such that 
\[
 d_\nu(w^{-1}u(g)\widehat{\gamma_g}(w),\id)<\varepsilon, \,\, \nu\in \Phi, \, g\in K.
\]
\end{lem}
\textbf{Proof.} Let $K\Subset G$, $\varepsilon>0$ be given. 
Take $\varepsilon'>0$ such that $8\sqrt{\varepsilon'}<\varepsilon$. Let $S'$ be a $(K,\varepsilon')$-invariant set, 
and set $S=S'\cup K$.
Choose $\delta>0$ such that $4|S'|\delta+8\sqrt{\varepsilon'}<\varepsilon$. 
Let  a cocycle crossed action $(\gamma,c)$, $ \Phi \Subset M_1(X,\mu)$,
and $u\in \Coc(G,[T])$ 
satisfying the condition
\[
 d_\nu(\widehat{\gamma_g}(u(s))^{-1}u(g)^{-1}u(gs), \id )<\delta, \,\, g,s\in S, \,\, \nu \in \Phi
\]
be given.
By Proposition \ref{prop:RohZschoice}, 
choose a partition $\{E_s\}_{s\in S'}$ of $X$, $E_s\supset Z_s$ and $w\in [T]$ such that 
\begin{align*}
 (1)\,\,\,& u(s)Z_s\subset E_s, \\
 (2)\,\,\,& \nu(\gamma_g(E_s\backslash Z_s))<\delta, \nu(E_s\backslash u(s)Z_s)<\delta,\,\, g\in K, \,\nu \in \Phi, \\
 (3)\,\,\,& \nu(\widehat{\gamma_g}(u(s))^{-1}u(g)^{-1}(E_{gs}\backslash Z_{gs}))<\delta,\,\, g\in K, \,s\in S'_g, \,\nu\in \Phi, \\
 (4)\,\,\,& \nu(E_{gs}\triangle \widehat{\gamma_g}(u(s))^{-1}g^{-1}E_{gs})<\delta,\,\,  g\in K, \,s\in S'_g,\,\nu\in \Phi, \\
(5)\,\,\,&  \sum_{s\in S'_g}\nu(E_{gs}\triangle \gamma_gE_s)<5\sqrt{\varepsilon'}, \,\,g\in K,\, \nu\in  \Phi, \\
 (6)\,\,\,&  \sum_{s\in S'\backslash S_{g^{-1}}'}\nu\left(E_s\right) <3\sqrt{\varepsilon'},\,\, g\in K, \,\nu \in \Phi, 
\\
 (7)\,\,\,& wx=u(s)x,\,\, x\in Z_s. 
\end{align*}
Let 
\[
 W_g:=\bigcup_{s\in S'_g}\{\widehat{\gamma_g}(u(s))^{-1}u(g)^{-1}u(gs)=\id\}\cap {\gamma_gZ_s}\cap\widehat{\gamma_g}(u(s))^{-1}u(g)^{-1}Z_{gs}.
\]
We can verify that $w^{-1}u(g)\widehat{\gamma_g}(w)=\id$ on $W_g$, $g\in K$,
 as in the proof of Lemma \ref{lem:appro2cochovanish1}. 

Next we show $\nu(X\backslash W_g)<\varepsilon$.
We have 
\[
 \nu(E_{gs}\triangle \gamma_gZ_s)\leq 
 \nu(E_{gs}\triangle \gamma_gE_s)+ \nu(\gamma_gE_{s}\backslash \gamma_gZ_s)<
 \nu(E_{gs}\triangle \gamma_gE_s)+ \delta
\]
by (2), and
\begin{align*}
\lefteqn{ \nu(E_{gs}\triangle \widehat{\gamma_g}(u(s))^{-1}u(g)^{-1}Z_{gs}) } \\
&\leq 
 \nu(E_{gs}\triangle \widehat{\gamma_g}(u(s))^{-1}u(g)^{-1}E_{gs})+
 \nu(\widehat{\gamma_g}(u(s))^{-1}u(g)^{-1}(E_{gs}\backslash  Z_{gs})) 
 <2\delta
 \end{align*}
by (3) and (4).
Hence we have
\[
  \nu\left(E_{gs}\triangle\left(\gamma_gZ_s\cap  \widehat{\gamma_g}(u(s))^{-1}u(g)^{-1}E_{gs}\right)\right)<3\delta +
\nu(E_{gs}\triangle \gamma_gE_s). 
\]
Then we have
\begin{align*}
\lefteqn{  \nu\left(\bigcup_{s \in  S'_g}E_{gs} \triangle
\bigcup_{s\in S'_g}\left(\gamma_gZ_s\cap  \widehat{\gamma_g}(u(s))^{-1}u(g)^{-1}E_{gs}\right)\right) } \\
&\leq \sum_{s\in S_g'}  \nu\left(E_{gs} \triangle
\left(\gamma_gZ_s\cap  \widehat{\gamma_g}(u(s))^{-1}u(g)^{-1}E_{gs}\right)\right) \\
&<\sum_{s\in S'_g}\left(3\delta+ \nu(E_{gs}\triangle \gamma_gE_s)\right) 
 <3|S'|\delta+5\sqrt{\varepsilon'}
\end{align*}
by (5).
Hence  we get
\begin{align*}
\lefteqn{ \nu\left(X\backslash \bigcup_{s\in S'_g}
\left(\gamma_gZ_s\cap  \widehat{\gamma_g}(u(s))^{-1}u(g)^{-1}E_{gs}\right)\right)} \\
&\leq\sum_{s\in S'\backslash S_{g^{-1}}'}\nu\left(E_s\right) +
\nu\left(\bigcup_{s\in S'_g}E_{gs} \triangle
\bigcup_{s\in S'_g}\left(\gamma_gZ_s\cap  \widehat{\gamma_g}(u(s))^{-1}u(g)^{-1}E_{gs}\right)\right) \\
 &< 3|S'|\delta+8\sqrt{\varepsilon'}
\end{align*}
by (6).
By the assumption
\[
d_\nu(\widehat{\gamma_g}(u(s))^{-1}u(g)^{-1}u(gs), \id )<\delta, \,\, g,s\in S, \,\nu\in \Phi,
\]
we have $\nu(X\backslash\bigcup_{s\in S'_g} \{\widehat{\gamma_g}(u(s))^{-1}u(g)^{-1}u(gs)=\id\})<|S'|\delta$.
Hence 
\[
 \nu(X\backslash W_g)<4|S'|\delta+8\sqrt{\varepsilon'}<\varepsilon
\]
holds. 
\hfill$\Box$

\begin{thm}\label{thm:2cohoVanish}
 Let $(\gamma,c)$ be an ultrafree cocycle crossed action of $G$. Then 
there exists $u\in \Coc(G,[T])$ such that 
${}_uc(g,h)=\id$, and hence
${}_u\gamma$ is an action. 

Moreover,
for any $e\in K\Subset G $, $\varepsilon>0$, there exists $S\Subset G$, $\delta>0$, which depends only on 
$K$ and $\varepsilon$, nor on cocycle crossed action $(\gamma,c)$, 
 such that
if 
\[
 d_\nu(c(g,h),\id)<\delta, \,\,g,h\in S, \,\nu\in \Phi
\]
for some $\Phi\Subset M_1(X,\mu)$ with $\mu\in \Phi$, 
then we can choose $u\in \Coc(G,[T])$ so that 
\[
 d_\nu(u(g),\id)<\varepsilon, \,\, g\in K, \,\, \nu\in \Phi.
\]
 \end{thm}
\textbf{Proof.} 
At first, we treat type II$_\infty$ or type III case.

Let  $e\in K\Subset G$ and $\varepsilon>0$ be given, and take $S\Subset G$ and $\delta>0$ 
as in Lemma \ref{lem:near1cocyclevanish}.
Assume 
$d_\nu(c(g,h),\id)<\delta$ for $g,h\in S$, $\nu \in \Phi\Subset M_1(X,\mu)$. 
There exists $v\in \Coc(G,[T])$ 
such that ${}_vc(g,h)=\id$ by Theorem \ref{2cohoinfinite}. 
Hence $c(g,h)=\widehat{\gamma_g}(v(h))^{-1}v(g)^{-1}v(gh)$ holds, and 
\[
 d_\nu(\widehat{\gamma_g}(v(h))^{-1}v(g)^{-1}v(gh),\id)<\delta,\,\, g,h\in S,\,\nu\in  \Phi.
\]
By Lemma \ref{lem:near1cocyclevanish}, there exists $w\in [T]$
such that 
\[
 d_\nu(w^{-1}v(g)\widehat{\gamma_g}(w),\id)<\varepsilon, \,\, \nu\in \Phi,\, g\in K.
\]
Define $u(g):=w^{-1}v(g)\widehat{\gamma_g}(w)$. Then 
we obtain
$d_\nu(u(g),\id)<\varepsilon$ for $g\in K$, $\nu \in \Phi$,
and 
\[
{}_uc(g,h)=
 u(g)\widehat{\gamma_g}(u(h))c(g,h)u(gh)^{-1}=
w^{-1} v(g)\widehat{\gamma_g}(v(h))c(g,h)v(gh)^{-1} w=\id.
\]
Hence we have proved the theorem for type II$_\infty$ and type III case.

Next, we assume $T$ is of type II$_1$. In this case, we can assume that 
$\mu$ is the unique $T$-invariant probability measure, and choose $\Phi$ as $\Phi=\{\mu\}$.
Let us take an increasing sequence $\{K_n\}_{n}\Subset G$, and decreasing sequence $\{\varepsilon_n\}_n$ such that
$e\in K_n$,  $\bigcup_{n=1}^\infty K_n=G$, and $\sum_{n}\varepsilon_n<\infty$. 
Take $S_n$ and $\delta_n$ for $K_n$ and $\varepsilon_n>0$  as in Lemma \ref{lem:near1cocyclevanish}.
We can choose $S_n$ and $\delta_n$ so that $S_n\subset S_{n+1}$, $\delta_n>\delta_{n+1}$.

For given $K\Subset G$, and $\varepsilon>0$, choose  $N\in \mathbb{N}$ 
such that $K\subset K_N$, $\varepsilon>\sum_{k=N}^\infty \varepsilon_k$. 
By Lemma \ref{lem:appro2cochovanish1}, 
take $S_N\Subset G$  and $\delta_N>0$ for $K_N$ and $\varepsilon_N>0$. 
Again by Lemma \ref{lem:appro2cochovanish1}, we can perturb $(\gamma,c)$ by some 
$w\in \Coc(G,[T])$
so that 
\begin{align*}
&d_\mu\left({}_wc(g,h),\id\right)<\varepsilon_N,\,\,g,h\in K_N, \,\,\, 
d_\mu\left({}_wc(g,h),\id\right)<\frac{\delta_N}{2},\,\, g,h\in S_N.
\end{align*}
Set 
\[
(\gamma^{(N)}, c_N):=({}_w\gamma,{}_wc), \,\,\,
u_{N}(g)=1.
\]

We will inductively construct a family of cocycle crossed actions $(\gamma^{(n)},c_n)$ 
and normalized maps $\{u_{n}\}\subset \Coc(G,[T])$
$n\geq N$, 
such that
\begin{align*}
(1.n)\,\,\,& (\gamma^{(n)},c_n)=({}_{u_{n}}\gamma^{(n-1)}, {}_{u_{n}}c_{n-1}),  \\
(2.n)\,\,\, & d_\mu(c_n(g,h),\id)<\varepsilon_n,\,\, g,h\in K_n, \\
(3.n)\,\,\,&  d_\mu(c_n(g,h),\id)<\frac{\delta_n}{2}, \,\,g,h \in S_n, \\
(4.n)\,\,\, & d_\mu(u_{n}(g),\id)<\varepsilon_{n-1},\,\, g\in K_{n-1}.\,\, 
\end{align*}
Here we regard $\gamma^{(N-1)}=\gamma^{(N)}$, $c_{N-1}(c,h)=c_N(g,h)$. 
Clearly  we have $(1.N)$, $(2.N)$, $(3.N)$ and $(4.N)$.

Assume we have done up to the $n$-th step.

By Lemma \ref{lem:appro2cochovanish1}, we choose $\bar{u}_{n+1}\in \Coc(G,[T])$
such that
\begin{align*}
(a.n+1)\,\,\, & d_\mu\left(\bar{u}_{n+1}(g)\widehat{\gamma^{(n)}_g}(\bar{u}_{n+1}(h))c_n(g,h)
  \bar{u}_{n+1}(gh)^{-1},\id\right)
<\varepsilon_{n+1}, \,\,
g,h\in K_{n+1}, \\
(b.n+1)\,\,\,  & d_\mu\left(\bar{u}_{n+1}(g)\widehat{\gamma^{(n)}_g}(\bar{u}_{n+1}(h))c_n(g,h)
\bar{u}_{n+1}(gh)^{-1},\id\right)
<\frac{\delta_{n+1}}{2},\,\, 
g,h\in S_{n+1}. 
\end{align*}
By $(b.n+1)$, we have
\[
d_\mu\left(\widehat{\gamma^{(n)}_g}(\bar{u}_{n+1}(h))^{-1}\bar{u}_{n+1}(g)^{-1}\bar{u}_{n+1}(gh), 
c_n(g,h)\right)<\frac{\delta_{n+1}}{2},\,\,\,g,h\in S_{n+1}.
\]
Combining with $(3.n)$, 
we get 
\[
d_\mu\left(\widehat{\gamma_g^{(n)}}(\bar{u}_{n+1}(h))^{-1}
\bar{u}_{n+1}(g)^{-1}\bar{u}_{n+1}(gh), \id\right)<\delta_n, \,\,\,
g,h\in S_n. 
\]
By Lemma \ref{lem:near1cocyclevanish}, there exists $w\in [T]$ such that 
$d_\mu(w^{-1}\bar{u}_{n+1}(g)\widehat{\gamma^{(n)}_g}(w),\id)<\varepsilon_n$ for $g\in K_n$.
Here set $u_{n+1}(g):=w^{-1}\bar{u}_{n+1}(g)\widehat{\gamma^{(n)}_g}(w)$. Then
we get $(4.n+1)$. Define a cocycle crossed action $(\gamma^{(n+1)}, c_{n+1})$ as $(1.n+1)$.
Then we get $(2.n+1)$ and $(3.n+1)$ from $(a.n+1)$ and $(b.n+1)$, respectively, 
and complete induction.

Let $v_{n}(g):=u_{n}(g)u_{n-1}(g)\cdots u_{N}(g)$. 
We have $(\gamma^n,c_n)=({}_{v_n}\gamma^{(N)}, {}_{v_n}c_N)$ by the construction, 
Fix $L\in \mathbb{N}$, and take any $g\in K_L$. 
By $(4.n)$, 
\[
 d_\mu(v_{n}(g),v_{n-1}(g))= d_\mu(u_{n}(g), \id)<\varepsilon_{n-1}, \,\,  n\geq L+1
\]
holds. So $\{v_n(g)\}_n$  is a Cauchy sequence, and hence 
$v_{n}(g)$ 
converges to some
$v(g)\in [T]$ 
uniformly. Note that $v_{n}(g)^{-1}$ converges to $v(g)^{-1}$ automatically, since 
$\mu$ is the invariant measure for $[T]$.
Combining with  $(2.n)$, we obtain ${}_vc(g,h)=\id$ for all $g,h\in G$.

If $g\in K_N$, then
\begin{align*}
 d_\mu(v_{n}(g),\id )&= d_\mu(v_{n}(g),v_{N}(g)) 
\leq   \sum_{k=N}^{n-1}d_\mu(v_{k+1}(g),v_{k}(g))
< \sum_{k=N}^{n-1}\varepsilon_k.
\end{align*}
Hence we have $d_\mu(v(g),\id)\leq \sum_{k=N}^\infty \varepsilon_k<\varepsilon$.
Set $S:=S_N\cup K_N$, $\delta:=\min\{\delta_N/2, \varepsilon_N\}$.
If $d_\mu(g,h)<\delta$ for $g,h\in S$, then we have $d_\mu(v(g),\id)<\varepsilon$ for $g\in K_N$.
Note that $S$ and $\delta$ are determined only on $K$ and $\varepsilon$.
\hfill$\Box$

\section{Classification}\label{sec:class}

\begin{lem}\label{lem:2cocyclecentral}
Let $\alpha$ and $\beta$ be actions of $G$ into $N[T]$ with 
$\md(\alpha_g)=\md(\beta_g)$. Then for any $\varepsilon>0$, $K\Subset G$, 
$\mu \in \Phi\Subset M_1(X,\mu)$, $\Lambda \Subset  [T]$, 
there exists $w\in \Coc(G,[T])$
such that \\
$(1)$ $\|{}_w\alpha_g(\xi)-\beta_{g}(\xi)\|<\varepsilon$, $g\in K$, $\xi\in \Phi$, \\
$(2)$ $d_\nu(\widehat{{}_w\alpha_g}(t), \widehat{\beta_g}(t))<\varepsilon$, $g\in K$, $t\in \Lambda,\nu\in\Phi $,  \\   
$(3)$ Let  $c(g,h):=w(g)\widehat{\alpha_g}(w(h))w(gh)^{-1}$. Then 
\[
 \|c(g,h)(\xi)-\xi\|<\varepsilon,\,\, d_\nu(\widehat{c(g,h)}(t),t)<\varepsilon,\,\, 
g,h\in K,\, \xi,\nu \in \Phi, \,t\in \Lambda.
\]
\end{lem}
\textbf{Proof.} 
By enlarging $K$, we may assume $e\in K=K^{-1}\Subset G$. 
Let 
\[
\tilde{\Phi}:=\left\{\beta_{gh}(\xi)\mid g,h\in K, \xi\in \Phi\right\}, \,\,\, 
\tilde{\Lambda}:=\left\{\widehat{\beta_{gh}}(t)\mid  g,h\in K,t\in \Lambda\right\}.
\]
By the assumption, $\beta_g\alpha_g^{-1}\in \Ker(\md)=\overline{[T]}$. Hence we can 
take $w\in \Coc(G,[T])$ so that 
\[
 \left\|{}_w\alpha_{gh}(\xi)-\beta_{gh}(\xi)\right\|<\frac{\varepsilon}{7},\,\,
d_\nu(\widehat{{}_w\alpha_{gh}} (t),\widehat{\beta_{gh}}(t))<\frac{\varepsilon}{7}
\]
for $g,h\in K$, $\nu, \xi\in \bigcup_{g\in K}\beta_g(\tilde{\Phi})$, 
$t\in \bigcup_{g\in K}\beta_g(\tilde{\Lambda})$.
Obviously, we have conditions (1), (2).

Then for $g,h\in K$, $\eta\in \tilde{\Phi}$, we have
\begin{align*}
 \left\|{}_w{\alpha}_{g}\,{}_w{\alpha}_{h}(\eta)-\beta_{gh}(\eta)\right\| &\leq 
 \left\|{}_w{\alpha}_{g}\,{}_w{\alpha}_{h}(\eta)-{_w}{\alpha}_{g}\beta_{h}(\eta)\right\| +
 \left\|{_w}{\alpha}_{g}\beta_{h}(\eta)-\beta_{g}\beta_{h}(\eta) \right\| \\
&\leq  \left\|{_w}{\alpha}_{h}(\eta)-\beta_{h}(\eta)\right\| +
 \left\|{_w}{\alpha}_{g}\beta_{h}(\eta)-\beta_{g}\beta_{h}(\eta) \right\| 
< \frac{2\varepsilon}{7}.
\end{align*}
Thus 
\begin{align*}
 \left\|c(g,h)\beta_{gh}(\eta)-\beta_{gh}(\eta)\right\| 
&\leq 
 \left\|c(g,h)\beta_{gh}(\eta)-c(g,h)\alpha_{gh}(\eta)\right\|+
 \left\|c(g,h)\alpha_{gh}(\eta)-\beta_{gh}(\eta)\right\| \\
&< \frac{3\varepsilon}{7}
\end{align*}
holds
for $g,h\in K$, $\eta\in \tilde{\Phi}$. Hence we get $\|c(g,h)(\xi)-\xi\|<3\varepsilon/7$ for $g,h\in K$, $\xi\in \Phi$.

For $g\in K$, $t\in \tilde{\Lambda}$, $\nu\in \tilde{\Phi}$, we have
\begin{align*}
 d_\nu(\widehat{c(g,h)}\widehat{{_w}{\alpha}_{gh}}(t), \widehat{\beta_{gh}}(t))  
&=d_\nu(\widehat{{_w}{\alpha}_g}\,\widehat{{_w}{\alpha}_h} (t) ,\widehat{\beta_{gh}}(t)) \\
&\leq
d_\nu(\widehat{{_w}{\alpha}_g}\,\widehat{{_w}{\alpha}_h} (t), 
\widehat{{_w}{\alpha}_g}\widehat{\beta_h}(t) )+
d_\nu(\widehat{{_w}{\alpha}_g}\widehat{\beta_h} (t), 
\widehat{\beta_{gh}}(t))
 \\
&\leq 
d_{{_w}{\alpha}_g^{-1}(\nu)}(\widehat{{_w}{\alpha}_h} (t),
\widehat{\beta_h}(t))+\frac{\varepsilon}{7}  \\
&\leq 
d_{\beta_g^{-1}(\nu)}(\widehat{{_w}{\alpha}_h }(t),
\widehat{\beta_h}(t) )+\|{_w}{\alpha}_g^{-1}(\nu)-\beta_g^{-1}(\nu)\|+
\frac{\varepsilon}{7} \\
&\leq 
d_{\beta_g^{-1}(\nu)}(\widehat{{_w}{\alpha}_h }(t),
\widehat{\beta_h}(t) )+\frac{2\varepsilon}{7} < \frac{3\varepsilon}{7}. 
\end{align*}

By noting $\|c(g,h)(\nu)-\nu\|\leq 3\varepsilon/7$ for $g,h\in K$, $\nu\in \Phi$, we have
\begin{align*}
 d_\nu( \widehat{c(g,h)}\widehat{\beta_{gh}}(t), \widehat{\beta_{gh}}(t))
&\leq
 d_\nu
(\widehat{c(g,h)}\widehat{\beta_{gh}}(t),
\widehat{ c(g,h)}\widehat{{_w}{\alpha}_{gh} }(t)) 
+
 d_\nu(\widehat{c(g,h)}\widehat{{_w}{\alpha}_{gh}} (t),\widehat{\beta_{gh}}(t)) \\
&\leq d_{c(g,h)^{-1}(\nu)}(\widehat{\beta_{gh}}(t),
\widehat{ {_w}{\alpha}_{gh}} (t)) 
+\frac{3\varepsilon}{7} \\
&\leq d_{\nu}(\widehat{\beta_{gh}}(t),
\widehat{ {_w}{\alpha}_{gh} }(t))
+\frac{6\varepsilon}{7} 
<\varepsilon
\end{align*}
for $\nu\in \Phi$, $g,h\in K$, $t\in \tilde{\Lambda}$.
Thus
$d_\nu(\widehat{c(g,h)}(t),t)<\varepsilon$ holds for $g,h\in K$, $t\in \Lambda$, $\nu\in \Phi$. \hfill$\Box$

\begin{lem}\label{lem:appro2cohovanish2}
 Let $\alpha$ and $\beta$ be actions of $G$ into $N[T]$ with $\md(\alpha_g)=\md(\beta_g)$. 
For any $\varepsilon>0$, $K\Subset G$, $\Lambda \Subset [T]$,
$\Phi\Subset  M_1(X,\mu)$, there exists $v\in \Coc(G,[T])$ such that
\begin{align*}
 & \|{_v}\alpha_g(\xi)-\beta_g(\xi)\|<\varepsilon,\,\, g\in K,\,\nu \in \Phi,  \\
&d_\nu(\widehat{{_v}\alpha_g}(t),\widehat{\beta_g}(t))<\varepsilon,\,\, g\in K,\,t\in \Lambda , \,\nu\in \Phi,\\
 &d_{\nu}(v(g)\widehat{\alpha_g}(v(h))v(gh)^{-1},\id)<\varepsilon,\,\, g,h \in K,\, \nu \in \Phi. 
\end{align*}
\end{lem}
\textbf{Proof.} 
Let $\tilde{\Phi}:=\{\beta_g(\xi)\mid g \in K, \xi\in \Phi\}$,
$\tilde{\Lambda }:=\{\widehat{\beta_g}(t)\mid g\in K, t\in \Lambda \}$.
Choose $\delta>0$ and $S$ for $\varepsilon/3>0$ and $K$ as in Lemma \ref{lem:appro2cochovanish1}.
By Lemma \ref{lem:2cocyclecentral}, there exists $u\in \Coc(G,[T])$
such that 
\begin{align*}
 &\|{}_u\alpha_g(\xi)-\beta_g(\xi)\|<\frac{\varepsilon}{3},\,\, g\in K,\,\xi \in \Phi, \\
& \|c(g,h)(\xi)-\xi\|<\delta,\,\, g,h \in S, \,\xi\in \tilde{\Phi},   \\
& d_\nu(\widehat{c(g,h)}(t),t)<\delta,\,\, g,h\in S,\,\, t\in \tilde{\Lambda },\, \nu\in \tilde{\Phi}, 
\end{align*}
where $c(g,h)=u(g)\alpha_g(u(h))u(gh)^{-1}$. 
By Lemma \ref{lem:appro2cochovanish1}, there exists $w\in\Coc(G, [T])$ such that
\[
 d_{\nu}(w(g)\widehat{{}_u{\alpha}_g}(w(h))c(g,h)w(gh)^{-1},\id)<\frac{\varepsilon}{3},\,\, g,h \in K, \,\nu \in \Phi
\]
and 
\[
  \|w(g)(\xi)-\xi\|<\frac{\varepsilon}{3},\,\,\, d_\nu(\widehat{w(g)}(t),t)<\frac{\varepsilon}{3},\,\, g\in K, \,
\xi,\nu\in \tilde{\Phi},\, t\in \tilde{\Lambda}.
\]
Let $v(g):=w(g)u(g)$. Then we have
\[
 d_{\nu}(v(g)\widehat{\alpha_g}(v(h))v(gh)^{-1},\id)<\varepsilon,\,\, g,h \in K,\, \nu \in \Phi. 
\]
 We can verify the first inequality as follows.
For $g\in K$,  $\xi\in \Phi$, 
\begin{align*}
 \|{}_v\alpha_g(\xi)-\beta_g(\xi)\|&\leq 
 \|w(g)u(g)\alpha_g(\xi)-w(g)\beta_g(\xi)\|+  \|w(g)\beta_g(\xi)-\beta_g(\xi)\| \\
&<\frac{2\varepsilon}{3}<\varepsilon
\end{align*} 
since $\beta_g(\xi)\in \tilde{\Phi}$.
Similarly, we have 
\begin{align*}
 d_\nu(\widehat{{_v}\alpha_g}(t),\widehat{\beta_g}(t))&\leq 
 d_\nu(\widehat{w(g){_u}\alpha_g}(t),\widehat{w(g)\beta_g}(t))+
 d_\nu(\widehat{w(g)\beta_g}(t),\widehat{\beta_g}(t)) \\
&\leq d_{w(g)(\nu)}(\widehat{{_u}\alpha_g}(t),\widehat{\beta_g}(t))
+\frac{\varepsilon}{3} \\
&\leq  \|w(g)(\nu)-\nu\|+
d_{\nu}(\widehat{{_u}\alpha_g}(t),\widehat{\beta_g}(t))  +\frac{\varepsilon}{3}
<\varepsilon
\end{align*}
for $g\in K$, $t\in \Lambda $, $\nu \in \Phi$. 
\hfill$\Box$

\begin{thm}\label{thm:approby1cocycle}
 Let $\alpha$ and $\beta$ be  ultrafree actions of $G$ into $N[T]$  with $\md(\alpha_g)=\md(\beta_g)$.
Then there exists a sequence $\{u_n(\cdot )\}$ of  1-cocycles for $\alpha_g$ such that 
$\lim\limits_{n\rightarrow \infty }{_{u_n}}\alpha_g=\beta_g$ in the $u$-topology.
\end{thm}
\textbf{Proof.} By Lemma \ref{lem:appro2cohovanish2},  
there exists  a sequence $\{v_n\} \subset  \Coc(G,[T])$ of 
normalized maps
such that 
$\lim\limits_{n\rightarrow \infty }{_{v_n}}\alpha_g=\beta_g$ in the $u$-topology, and 
$\lim\limits_{n\rightarrow \infty}d_\mu\left(v_n(g)\widehat{\alpha_g}(v_n(h))v_{n}(gh)^{-1},\id\right)=0$.
Let $\alpha^{(n)}={_{v_n}}\alpha$ and $c_n(c,h)=v_n(g)\widehat{\alpha_g}(v_n(h))v_{n}(gh)^{-1}$.
By Theorem \ref{thm:2cohoVanish}, there exists a sequence
$\{w_n\} \subset  \Coc(G,[T])$
such that 
\[ 
w_n(g)\widehat{\alpha_g^{(n)}}(w_n(h))c_n(g,h)w_{n}(gh)^{-1}=1, \,\,
\lim_{n\rightarrow \infty} d_\mu(w_n(g),\id)=0.
\]
Then it turns out that $u_n(g):=w_n(g) v_n(g)$ is a 1-cocycle for $\alpha_g$, 
and $\lim\limits_{n\rightarrow \infty}{_{u_n}}\alpha_g=\beta_g$ holds in the $u$-topology. \hfill$\Box$

\begin{lem}\label{lem:1cohovanish}
 Let $K\Subset G$ and $\varepsilon>0$ be given. 
Then there exist $S\Subset G$ and $\delta>0$ satisfying the following; 
for any action $\gamma$ of $G$, a 1-cocycle $u(\cdot)$ for $\gamma$, $\Phi \Subset M_1(X,\mu)$ with $\mu\in \Phi$ 
and $\Lambda  \Subset [T]$ 
satisfying 
\[\|u(s)(\xi)-\xi\|<\delta,\,\,\,  d_\nu(\widehat{u(s)}(t),t)<\delta,\,\,\,
s\in S, \,\, \xi, \nu\in \Phi, \,\, t\in \Lambda,  
\]
there exists $ w\in [T]$ such that 
\[
 d_\nu(u(g)\widehat{\gamma_g}(w)w^{-1},1)<\varepsilon, 
\,\,\, 
 \|w(\xi)-\xi\|<\varepsilon, \,\,\, d_\nu(w(t),t)<\varepsilon,\,\,\, g\in K, \,\, \xi,\nu\in \Phi,\,\, t\in \Lambda .
\]
\end{lem}
\textbf{Proof.} Take $\varepsilon_1>0$ with $8\varepsilon_1^{\frac{1}{2}}<\varepsilon$, and 
let $S$ be a $(K,\varepsilon_1)$-invariant set.
Choose $\delta>0$ with $8\varepsilon_1^\frac{1}{2}+3|S|\delta<\varepsilon$, 
$4|S|\delta<\varepsilon$. 

By Proposition \ref{prop:RohZschoice}, 
take a partition $\{E_s\}_{s\in S} $ of $X$, $Z_s\subset E_s$ and $w\in [T]$
such that 
\begin{align*}
(1)\,\,\, & u(s) Z_s\subset E_s, \,\, s\in S,  \\
(2)\,\,\, & \nu(E_s\backslash Z_s)<\delta, \nu(E_s\backslash u(s) Z_s)<\delta, \,\, s\in S, \,\nu\in \Phi,\\
(3)\,\,\,& \nu\left(u(gs)\gamma_g\left(E_{s}\backslash Z_s\right)\right)<\delta,\,\, g\in K,  
       \,s\in S_g, \, \nu\in \Phi,\\
(4)\,\,\,& \nu(u(s)t^{-1}(E_s\backslash  Z_s))<\delta, \,\, s\in S,\, t\in \Lambda, \,\nu\in \Phi,\\
(5)\,\,\,& \nu(u(s)E_s\triangle E_s)<\delta, \,\, s\in S, \,\nu\in \Phi,\\
(6)\,\,\,&\nu(E_{gs}\triangle \widehat{\gamma_g} (u(s))E_{gs}) <\delta,\,\, s\in S_g,\, \nu\in \Phi,\\
(7)\,\,\,& \nu(E_s\triangle u(s) t^{-1}E_s)<\delta, \,\, s\in S, \,t\in \Lambda,\, \nu\in \Phi,\\
(8)\,\,\,& \sum_{s\in S_g}u(gs)^{-1}(\nu)(\gamma_gE_s\triangle E_{gs} )<5\varepsilon_1^{\frac{1}{2}}, \,\, 
    g\in K, \,\nu\in \Phi,\\
(9)\,\,\,&\sum_{s\in S\backslash S_{g^{-1}}}\nu(E_s)<3\varepsilon_1^{\frac{1}{2}}, \,\, g\in K,\\  
(10)\,\,\, & wx=u(s)x,\,\, x\in Z_s. 
\end{align*}

In the following proof, the letter $g$, $s$, and $\nu$ 
denote an element in $K$, $S$, and $\Phi$, respectively.
As in the proof of Lemma \ref{lem:appro2cochovanish1}, 
we can see that  
\[
 u(g)\widehat{\gamma_g}(w)w^{-1} x=u(g)\gamma_gu(s)\gamma_g^{-1}u(gs)^{-1}x=x
\]
for $x\in u(gs)Z_{gs}\cap u(gs)\gamma_gZ_s$.

We have
\begin{align*}
\lefteqn{ \nu\left(E_{gs}\triangle \left(u(gs)Z_{gs}\cap u(gs)\gamma_gZ_s\right)\right)} \\
&\leq 
\nu\left(E_{gs}\backslash  u(gs)Z_{gs}\right)+
\nu\left(E_{gs}\triangle u(gs)\gamma_gZ_s\right) \\
&< \delta+ 
\nu\left(E_{gs}\triangle u(gs)\gamma_g E_s\right) +
\nu\left(u(gs)\gamma_g(E_{s})\backslash  u(gs)\gamma_gZ_s\right)\,\,\, (\mbox{by (2)}) \\
&< 2\delta +
\nu\left(E_{gs}\triangle u(gs) E_{gs}\right) +
\nu\left(u(gs)\left(E_{gs}\triangle \gamma_g E_s\right)\right)\,\,\, (\mbox{by (3)}) \\
 &<
3\delta +
u(gs)^{-1}(\nu)\left(E_{gs}\triangle \gamma_g E_s\right) \,\,\, (\mbox{by (5)}).
\end{align*}
Thus 
\begin{align*}
\lefteqn{\nu\left(X\backslash \bigcup_{s\in S_g}u(gs)Z_{gs}\cap u(gs)\gamma_gZ_s\right) } \\
&\leq \nu(X\backslash \bigsqcup_{s\in S_g}E_{gs})+
\sum_{s\in S_g}\nu\left(E_{gs}\triangle \left(u(gs)Z_{gs}\cap u(gs)\gamma_gZ_s\right)\right) \\
&\leq \sum_{s \in S\backslash S_{g^{-1}}}\nu(E_s)+
\sum_{s\in S_g}\left(3\delta +
u(gs)^{-1}(\nu)\left(E_{gs}\triangle \gamma_g E_s\right) \right)\\
&< 3\varepsilon_1^{\frac{1}{2}}+3|S|\delta +5\varepsilon_1^{\frac{1}{2}} 
= 8\varepsilon_1^{\frac{1}{2}}+3|S|\delta 
<\varepsilon
\end{align*}
holds.
Hence 
$\nu\left(\left\{u(g)\widehat{\gamma_g} (w)w^{-1}\ne \id\right\}\right)<\varepsilon$ for 
$g\in K$ and $\nu\in \Phi$, which implies 
\[
 d_\nu(u(g) \widehat{\gamma_g} (w)w^{-1},\id)<\varepsilon,\,\, g\in K,\,\, \nu \in \Phi.
\]

We next show $\|w(\xi)-\xi\|<\varepsilon$ and $d_\nu(\widehat{w}(t),t)<\varepsilon$.
Let $Z=\bigsqcup_{s\in S}Z_s$.
As in the proof of Lemma \ref{lem:appro2cochovanish1},
we can see $w(\xi)(x)=u(s)(\xi)(x)$ on $u(s)Z_s$, 
and 
\[
 \int_{X\backslash wZ}|w(\xi)(x)-\xi(x)|d \mu(x) <2|S|\delta
\]
by using (2) and (10).
If $u(s)$ satisfies $\|u(s)(\xi) -\xi \|<\delta$ for $s\in S$, then
\begin{align*}
 \|w(\xi) -\xi \|
&=\sum_{s\in S}\int_{u(s)Z_s}|w(\xi)(x)-\xi(x)|d \mu(x) +\int_{X\backslash wZ}|w(\xi)(x)-\xi(x)|d \mu(x)  \\
&<  \sum_{s\in S} \int_{u(s)Z_s}|u(s)(\xi)(x)-\xi(x)|d \mu(x) +2|S|\delta 
< 3|S|\delta 
 <\varepsilon
\end{align*}
holds for $\xi\in \Phi$.

For $t\in \Lambda \subset [T]$, and $x\in u(s)Z_s\cap u(s) t^{-1}Z_s$, 
$w^{-1}x=u(s)^{-1}x\in Z_s\cap  t^{-1}Z_s$. 
Hence 
$tw^{-1}x=u(s)^{-1}x\in tZ_s\cap  Z_s$, and $wtw^{-1}x=u(s)tu(s)^{-1}x$ holds. 

Then 
\begin{align*}
 \nu\left(E_s\triangle \left(u(s)Z_s\cap u(s) t^{-1}Z_s\right)\right)&\leq 
 \nu(E_s\backslash  u(s)Z_s)+
 \nu(E_s\triangle u(s) t^{-1}Z_s) \\
&< \delta + \nu(E_s\triangle u(s) t^{-1}E_s)+ \nu(u(s)t^{-1}(E_s\backslash  Z_s)) \,\,\, (\mbox{by (2)}) \\
 &< 3\delta \,\,\, (\mbox{by (4) and (7)}). \\ 
\end{align*}

 Let us assume $d_\nu(\widehat{u(s)}(t),  t )<\delta$. Hence 
$A_{s,t}:=\{\widehat{u(s)}(t)=  t \}$ satisfies $\nu(X\backslash A_{s,t})<\delta$.
Thus 
\begin{align*}
\lefteqn{\nu\left(X\backslash 
\bigcup_{s\in S}\left(u(s)Z_s\cap u(s) t^{-1}Z_s\cap A_{s,t}\right)\right)} \\
&\leq \sum_{s\in S}  \nu\left(E_s\triangle \left(u(s)Z_s\cap u(s) t^{-1}Z_s\right)\right)
+\sum_{s\in S}\nu(X\backslash A_{s,t}) \\
&\leq 4\delta|S| <\varepsilon
\end{align*}
and we have
$\nu(\{\widehat{w}(t) \ne  t \})<\varepsilon$, equivalently $d_\nu(\widehat{w}(t),t)<\varepsilon$.
\hfill$\Box$

\medskip

\noindent
\textbf{Remark.} In Lemma \ref{lem:1cohovanish}, we can choose $\delta$ and $S$ so that 
$\delta<\delta'$ and $S'\subset S$ for any given $\delta'>0$ and $S'\Subset G$.

\smallskip

Now we can classify ultrafree actions.
\begin{thm}\label{thm:classfree}
 Let $\alpha$ and $\beta$ be ultrafree actions of $G$ into $N[T]$ with $\md(\alpha_g)=\md(\beta_g)$. Then they 
are strongly cocycle conjugate.
\end{thm}
\noindent
\textbf{Proof.}
Let $\{\xi_i\}_{i=0}^\infty$ be a countable dense subset of $M_1(X,\mu)$ with $\xi_0=\mu$.
Take $\varepsilon_n>0$ and $K_n\Subset G$ 
 such that $\sum_{n=0}^\infty \varepsilon_n<\infty$, $\varepsilon_n>\varepsilon_{n+1}, $
$e\in K_{n}$, $K_n\subset K_{n+1}$,  $\bigcup_{n=0}^\infty K_n=G$.
Then choose $S_n\Subset G$, $\delta_n>0$ for $K_n$, $\varepsilon_n$  as in  Lemma \ref{lem:1cohovanish}. 
We can assume $S_n\subset S_{n+1}$ and $\delta_{n+1}<\delta_n$. (See a remark after Lemma \ref{lem:1cohovanish}.)

Set $\gamma_g^{(0)}:=\alpha_g$, $\gamma^{(-1)}_g:=\beta_g$, and 
construct actions $\gamma_g^{(n)}$ of $G$, 
$v_{n}(g), \bar{v}_{n}(g), w_n,\theta_n \in [T]$, 
$\Phi_n \Subset M_1(X,\mu)$ and $\Lambda_n\Subset [T]$ as follows;
\begin{align*}
(1.n) \,\, & \gamma_g^{(n)}=\bar{v}_{n}(g)w_n\gamma_g^{(n-2)}w_n^{-1}, \\
(2.n) \,\, & \theta_n=w_n\theta_{n-2}, \\
(3.n)\,\, & v_{n}(g)= \bar{v}_{n}(g)\widehat{w_n}(v_{n-2}(g)),  \\
(4.n)\,\, & \|\gamma_g^{(n)}(\xi)-\gamma_g^{(n-1)}(\xi)\|<\varepsilon_n,\,\, g\in K_n,\,\xi\in \Phi_{n-1}, \\
 (5.n)\,\,  & d_\mu\!\left(\widehat{\gamma_g^{(n)}}(t),\widehat{\gamma_{g}^{(n-1)}}(t)\right)<\varepsilon_n, 
\,\, g\in K_n,\, t\in \Lambda_{n-1},\\
 (6.n)\,\,  & \|\gamma_g^{(n)}(\xi)-\gamma_g^{(n-1)}(\xi)\|<\frac{\delta_{n-1}}{2}, \,\,
g\in S_{n-1}, \,  \xi\in \bigcup_{g\in S_{n-1}}\gamma^{(n-1)}_{g^{-1}}(\Phi_{n-1}), \\
 (7.n)\,\, & 
d_\nu\!\left(\widehat{\gamma_g^{(n)}}(t),\widehat{\gamma_g^{(n-1)}}(t)\right)<\frac{\delta_{n-1}}{2}, \,\,
      g\in S_{n-1}, \,
t \in \bigcup_{s\in S_{n-1}}\gamma^{(n-1)}_{g^{-1}}(\Lambda_{n-1}), \, 
\nu  \in \Phi_{n-1}, \\
 (8.n)\,\,& d_\nu(\bar{v}_{n}(g),\id)<\varepsilon_{n-2},\,\, g\in K_{n-2},\, \nu\in \Phi_{n-2}, 
\,\,\, (n\geq 2),\\
 (9.n)\,\,& \|w_n(\xi)-\xi\|<\varepsilon_{n-2},\,\, \xi \in \Phi_{n-2},\,\,\, (n\geq 2),\\ 
(10.n)\,\, &d_\nu(\widehat{w_n}(t),t)<\varepsilon_{n-2},\,\, \nu\in \Phi_{n-2},\, t\in \Lambda_{n-2}, \,\,\,(n\geq 2), \\ 
(11.n)\,\, & \Phi_n=\{\xi_i\}_{i=0}^n \cup \{\theta_{n}(\xi_i)\}_{i=0}^n\cup \{v_{n}(g)(\mu)\}_{g\in K_n}, \\ 
 (12.n)\,\, & \Lambda_n=\{T^i\}_{i=-n}^n \cup \{\theta_n(T^i)\}_{i=-n}^n \cup \left\{v_{n}(g), 
v_{n}(g)^{-1}\right\}_{g\in K_n}. 
\end{align*}

\noindent
1st step.  
Let $\theta_{-1}=\theta_0=\id$, 
$v_{-1}(g)=v_{0}(g)=\id$.
By Theorem \ref{thm:approby1cocycle}, take a 1-cocycle $u_{1}(\cdot )$ for $\gamma^{(-1)}$such that 
\begin{align*}
 (a.1)\,\, & 
\|{_{u_1}}\gamma_g^{(-1)}(\xi)-\gamma_g^{(0)}(\xi)\|<\varepsilon_{1},\,\,  g\in K_{1}, \,
\xi \in \Phi_{0}, \\
 (b. 1)\,\,&  
d_\mu \!\left(\widehat{{_{u_1}}\gamma_g^{(-1)}}(t), \widehat{\gamma_g^{(0)}}(t)\right)<\varepsilon_{1}, \,\, g\in K_{1}, \,
 t\in \Lambda_{0},  \\
 (c.1)\,\, &
\|{_{u_1}}\gamma_g^{(-1)}(\xi)-\gamma_g^{(0)}(\xi)\|<\frac{\delta_{0}}{2}, \,\,g\in S_{0}, \,
\xi \in \bigcup_{g\in S_{0}}\gamma_{g^{-1}}^{(0)}(\Phi_{0}), \\
(d. 1) \,\, &d_\nu\!\left(\widehat{{_{u_1}}\gamma_g^{(-1)}}(t), \widehat{\gamma_g^{(0)}}(t)\right)<\frac{\delta_{0}}{2},   \,\,
     g\in S_{0}, \,
t \in \bigcup_{g\in S_{0}}\gamma^{(0)}_{g^{-1}}(\Lambda_0),  \,
\nu \in \Phi_{0}.
\end{align*} 
Set $w_1=\id$, $\bar{v}_{1}(g)=u_{1}(g)$, and 
define
\begin{align*}
\gamma_g^{(1)}& :=\bar{v}_{1}(g)w_1\gamma_g^{(-1)}w_1^{-1}   
={_{u_1}}\gamma_g^{(-1)}, \\
  \theta_1&:=w_1\theta_{-1}=\id,  \\
v_{1}(g)&:=\bar{v}_{1}(g)\widehat{w_1}(v_{-1}(g))=u_{1}(g) 
\end{align*} 
as in $(1.1)$, $(1.2)$, $(1.3)$, respectively.
By $(a.1)$, $(b.1)$, $(c.1)$, $(d.1)$, we get $(4.1)$, $(5.1)$, $(6.1)$ and $(7.1)$, respectively.
Define $\Phi_1$ and $\Lambda_1$ as in $(11.1)$, $(12.1)$, respectively. 
Then we finished the 1st step of induction.

Assume that we have done up to the $n$-th step.
By Theorem \ref{thm:approby1cocycle}, let us take a $\gamma^{(n-1)}$-cocycle $u_{n+1}(\cdot)$ such that 
\begin{align*}
(a.n+1)\,\, & 
\|{_{u_{n+1}}}\gamma_g^{(n-1)}(\xi)-\gamma_g^{(n)}(\xi)\|<\varepsilon_{n+1},\,\,  g\in K_{n+1}, \,
\xi \in \Phi_{n}, \\
 (b. n+1)\,\,&  
d_\mu\! \left(\widehat{{_{u_{n+1}}}\gamma_g^{(n-1)}}(t), \widehat{\gamma_g^{(n)}}(t)\right)
<\varepsilon_{n+1}, \,\, g\in K_{n+1}, 
 t\in \Lambda_{n},  \\
 (c.n+1)\,\, &
\|{_{u_{n+1}}}\gamma_g^{(n-1)}(\xi)-\gamma_g^{(n)}(\xi)\|<\frac{\delta_{n}}{2}, \,\,g\in S_{n}, \,
\xi \in \bigcup_{g\in S_{n}}\gamma_{g^{-1}}^{(n)}(\Phi_{n}), \\
(d. n+1) \,\, &d_\nu\!\left(
\widehat{{_{u_{n+1}}}\gamma_g^{(n-1)}}(t), \widehat{\gamma_g^{(n)}}(t)\right)<\frac{\delta_{n}}{2},   \,\,
     g\in S_{n}, \,
t \in \bigcup_{g\in S_{n}}\gamma^{(n)}_{g^{-1}}(\Lambda_n),  \,
\nu \in \Phi_{n}, \\
 (e.n+1)\,\, & 
\|{_{u_{n+1}}}\gamma_g^{(n-1)}(\xi)-\gamma_g^{(n)}(\xi)\|<\frac{\delta_{n-1}}{2},\,\, g\in S_{n-1}, \,
\xi \in \bigcup_{g\in S_{n-1}}\gamma_{g^{-1}}^{(n-1)}(\Phi_{n-1}), \\
(f. n+1) \,\, &d_\nu\!\left(\widehat{{}_{u_{n+1}}\gamma_g^{(n-1)}}(t), 
\widehat{\gamma_g^{(n)}}(t)\right)<\frac{\delta_{n-1}}{2},   \\ \,\,
      & \hspace*{4.8cm} g\in S_{n-1}, \, t \in \bigcup_{g\in S_{n-1}}\widehat{\gamma^{(n-1)}_{g^{-1}}}(\Lambda_{n-1}),  \,
\nu \in \Phi_{n-1}. 
\end{align*}

By $(6.n)$ and $(e.n+1)$, we have
\[
 \|u_{n+1}(g)\gamma_g^{(n-1)}(\xi)-\gamma_g^{(n-1)}(\xi)\|<\delta_{n-1},\, \,
g\in S_{n-1}, \,   \xi\in \bigcup_{g\in S_{n-1}}\gamma^{(n-1)}_{g^{-1}}(\Phi_{n-1})
\]
and hence
\[
 \|u_{n+1}(g)(\xi)-\xi\|<\delta_{n-1},\,\, 
g\in S_{n-1},\,   \xi\in \Phi_{n-1}.
\]

By $(7.n)$ and $(f.n+1)$,
\[
 d_\nu\!\left(\widehat{u_{n+1}(g)}\widehat{\gamma_g^{(n-1)}}(t),\widehat{\gamma_g^{(n-1)}}(t)\right)<\delta_{n-1}, \,\,
      g\in S_{n-1}, \,
t \in \bigcup_{g\in S_{n-1}}\widehat{\gamma^{(n-1)}_{g^{-1}}}(\Lambda_{n-1}),  \,
\nu \in \Phi_{n-1},  
\]
and hence
\[
 d_\nu\!\left(\widehat{u_{n+1}(g)}(t),t\right)<\delta_{n-1}, \,\, g\in S_{n-1}, \,
t \in \Lambda_{n-1}, \,
\nu \in \Phi_{n-1}.
\]

By Lemma \ref{lem:1cohovanish}, there exists $w_{n+1}\in [T]$ such that
\[
 d_\nu\!\left(u_{n+1}(g)\widehat{\gamma^{(n-1)}_g}(w_{n+1})w_{n+1}^{-1},\id\right)<\varepsilon_{n-1},   
\,\, g\in K_{n-1}, \, \nu\in \Phi_{n-1},
\]
\[
 \|w_{n+1}(\xi)-\xi\|<\varepsilon_{n-1}, \,\,\, 
d_\nu(\widehat{w_{n+1}}(t),t)<\varepsilon_{n-1}, \,\, \xi,\nu\in \Phi_{n-1}, 
\,t\in K_{n-1}.
\]

Set 
\begin{align*}
\bar{v}_{n+1}(g)&:=u_{n+1}(g)\widehat{\gamma^{(n-1)}_g}(w_{n+1})w_{n+1}^{-1},  \\
 \gamma^{(n+1)}_g&:={_{u_{n+1}}}\gamma_g^{(n-1)} =\bar{v}_{n+1}(g)w_{n+1}\gamma^{(n-1)}_g w_{n+1}^{-1}, \\
 \theta_{n+1}&:=w_{n+1}\theta_{n-1}. 
 \end{align*}
We clearly have $(1.n+1)$, $(2.n+1)$, $(3.n+1)$, $(8.n+1)$, $(9.n+1)$, and $(10.n+1)$.
From  $(a.n+1)$, $(b.n+1)$, $(c.n+1)$ and $(d.n+1)$, 
we obtain $(4.n+1)$,  $(5.n+1)$, $(6.n+1)$ and $(7.n+1)$, respectively.
We define $\Phi_{n+1}$ and $\Lambda_{n+1}$ as in $(11.n+1)$ and $(12.n+1)$, respectively.
Then we finished the $(n+1)$-st step, and completed induction. 

By the construction, we have
\[
\gamma_{g}^{(2n)}=v_{2n}(g)\theta_{2n}\alpha_{g}\theta_{2n}^{-1},\,\,\, 
\gamma_{g}^{(2n+1)}=
v_{2n+1}(g)\theta_{2n+1}\beta_{g}\theta_{2n+1}^{-1}.
\]
We will show that sequences $\{\theta_{2n}\}_n$, $\{\theta_{2n+1}\}_n$, $\{v_{2n}(g)\}_n$ and 
$\{v_{2n+1}(g)\}_n$ will converge.
Fix $k\in \mathbb{N}$, and take  $\xi\in \{\xi_i\}_{i=1}^k$, $t\in \{T^l\}_{|l|\leq k}$.
 For  $n>k+2$, we have $\xi, \theta_{n-2}(\xi)\in \Phi_{n-2}$, $ \widehat{\theta_{n-2}}(t)\in \Lambda_{n-2}$.
Then
\[
 \|\theta_{n}(\xi)-\theta_{n-2}(\xi)\|=\|w_n(\theta_{n-2}(\xi))-\theta_{n-2}(\xi)\|<\varepsilon_{n-2},
\]
\[
 \|\theta_{n}^{-1}(\xi)-\theta_{n-2}^{-1}(\xi)\|=\|w_n^{-1}(\xi)-\xi\|<\varepsilon_{n-2}
\]
and 
\[
 d_\mu\!\left(\widehat{\theta_n}(t),\widehat{\theta_{n-2}}(t)\right)=
d_\mu\!\left(\widehat{w_n}(\widehat{\theta_{n-2}}(t)),\widehat{\theta_{n-2}}(t)\right)<\varepsilon_{n-2}
\]
hold by $(9.n)$ and $(10.n)$.
It follows that $\{\theta_{2n}\}_n$ and $\{\theta_{2n+1}\}_n$ are both
Cauchy sequences with respect to the metric $d$ on $N[T]$.
(See \S \ref{subsec:topology} on the definition of $d$.)
Hence both $\{\theta_{2n}\}_n$ and $\{\theta_{2n+1}\}_n$ converge to some $\sigma_0,\sigma_1\in \overline{[T]}$, 
respectively in the $u$-topology. 

Fix $l\in \mathbb{N}$ and take any $g\in K_l$. Then for  $n>l+2$, we have $v_{n-2}(g),  v_{n-2}(g)^{-1}
\in \Lambda_{n-2}$, 
$v_{n-2}(g)(\mu)\in \Phi_{n-2}$. Thus
\begin{align*}
\lefteqn{ d_\mu(v_{n}(g),v_{n-2}(g))} \\
&\leq  d_\mu(\bar{v}_{n}(g)\widehat{w_{n}}(v_{n-2}(g)), 
\bar{v}_{n}(g)v_{n-2}(g))+ d_\mu(\bar{v}_{n}(g)v_{n-2}(g),v_{n-2}(g))\\ 
&=  d_\mu\left(\widehat{w_{n}}(v_{n-2}(g)), 
v_{n-2}(g)\right)+ d_{v_{n-2}(g)(\mu)}(\bar{v}_{n}(g),\id)\\
&< 2\varepsilon_{n-2}
\end{align*}
and 
\begin{align*}
\lefteqn{d_\mu\!\left(v_{n}(g)^{-1},v_{n-2}(g)^{-1}\right)} \\
&\leq  d_\mu\!\left(\widehat{w_n}\left(v_{n-2}(g)^{-1}\right)\bar{v}_{n}(g)^{-1},
\widehat{w_n}\left(v_{n-2}(g)^{-1}\right)\right) 
+
d_\mu\!\left(\widehat{w_n}\left(v_{n-2}(g)^{-1}\right),v_{n-2}(g)^{-1}\right)\\ 
&=  d_\mu\left(\bar{v}_{n}(g)^{-1}, \id\right)
+d_\mu\!\left(\widehat{w_n}\left(v_{n-2}(g)^{-1}\right),v_{n-2}(g)^{-1}\right)\\ 
&<2\varepsilon_{n-2}
\end{align*}
by $(8.n)$ and $(10.n)$.
Thus both $\{v_{2n}(g)\}_n$ and $\{v_{2n+1}(g)\}_n$ are Cauchy sequences with respect to $d_\mu$, 
and hence 
converge to some $z_0(g), z_1(g)\in [T]$
uniformly, respectively. 

Summarizing these results, we have
 \begin{align*}
\lim_{n\rightarrow \infty}\gamma_{g}^{(2n)}&=
\lim_{n\rightarrow \infty}v_{2n}(g)\theta_{2n}\alpha_{g}\theta_{2n}^{-1}=
z_0(g)\sigma_0\alpha_g\sigma_0^{-1}, \\
\lim_{n\rightarrow \infty}\gamma_{g}^{(2n+1)}&=
\lim_{n\rightarrow \infty}v_{2n+1}(g)\theta_{2n+1}\beta_{g}\theta_{2n+1}^{-1}=
z_1(g)\sigma_1\beta_g\sigma_1^{-1}.
 \end{align*}

By $(4.n)$ and $(5.n)$, 
we have $z_0(g)\sigma_0\alpha_g\sigma_0^{-1}=z_1(g)\sigma_1\beta_g\sigma_1^{-1}$.
Hence $\alpha$ and $\beta$ are cocycle conjugate. \hfill$\Box$

\medskip

\noindent
\textbf{Proof of Theorem \ref{thm:main}.}
Let $N:=N_\alpha=N_\beta$, $Q:=G/N$, and $\pi:G\rightarrow Q$ be the quotient map.
Fix a section $s:Q\rightarrow G$ such that $s(e)=e$.
Then $\alpha_{s(p)}$ is an ultrafree cocycle crossed action of $Q$.
By Theorem \ref{thm:2cohoVanish}, there exists  $v\in \Coc(Q,[T])$
such that  $\bar{\alpha}_p:=v(p)\alpha_{s(p)}$ is a genuine action of $Q$. 
Here define $v(g):=v(p)\alpha_n^{-1}\in [T]$, where $g=ns(p)$ with $p=\pi(g)$ and  $n\in N$. 
Then $v(g) \alpha_g=v(p)\alpha_{s(p)}=\bar{\alpha}_{\pi(g)}$, and $\bar{\alpha}_{\pi(g)}$ 
is an action of $G$. 
Thus $\alpha_g$ is strongly cocycle conjugate to $\bar{\alpha}_{\pi(g)}$ for some ultrafree action $\bar{\alpha}$
of $Q$. In the same way, $\beta_g$ 
is strongly cocycle conjugate to $\bar{\beta}_{\pi(g)}$ for some ultrafree action $\bar{\beta}$
of $Q$. Since $\md(\bar{\alpha}_p)=\md(\bar{\beta}_p)$, $\bar{\alpha}$ and 
$\bar{\beta}$ are strongly cocycle conjugate as actions of $Q$ by Theorem \ref{thm:classfree},
and hence so are as actions of $G$. 
Therefore two actions $\alpha$ and $\beta$ of $G$ are strongly cocycle conjugate. \hfill$\Box$

\ifx\undefined\bysame
\newcommand{\bysame}{\leavevmode\hbox to3em{\hrulefill}\,}
\fi

\end{document}